\documentclass[a4paper,12pt]{article}
\usepackage[applemac]{inputenc}
\usepackage{amssymb}
\usepackage{graphicx}
\usepackage{amsmath}
\usepackage{calc}
\usepackage{eurosym}
\usepackage[english]{babel}
\usepackage{calc} 
\begin{document}

\title{\textbf{A phase transition for the heights of a fragmentation tree}}

\author{\textbf{Adrien Joseph}}

\date{}

\maketitle

\begin{center}
\emph{Laboratoire de Probabilités et Modèles Aléatoires \ \\
Université Pierre et Marie Curie \ \\
175, rue du Chevaleret \ \\
F-75013 Paris, France}
\end{center}

\newtheorem{lemme}{Lemma}
\newtheorem{theoreme}{Theorem}
\newtheorem{prop}{Proposition}
\newtheorem{coro}{Corollary}
\newtheorem{defi}{Definition}
\newtheorem{remark}{Remark}

\ \\
\textbf{Summary.} We provide information about the asymptotic regimes for a homogeneous fragmentation of a finite set. We establish a phase transition for the asymptotic behaviours of the shattering times, defined as the first instants when all the blocks of the partition process have cardinality less than a fixed integer. Our results may be applied to the study of certain random split trees.
 

\ \\
\textbf{A.M.S. Classiﬁcation.} 60 F 15, 60 J 80.

\ \\
\textbf{e-mail.} \verb+adrien.joseph@etu.upmc.fr+

\section{Introduction}

The occupancy scheme is a simple urn model which often occurs in probability, statistics, combinatorics and computer science. We can cite for example species sampling \cite{bunfitz, hegast}, analysis of algorithms \cite{gardy}, learning theory \cite{bouchegar}, etc. The books by Johnson and Kotz \cite{johsonkotz} and by Kolchin \emph{et al.}\,\,\cite{kolchinat} are standard references. The present work is partly motivated by the study of fragmentation trees.

Let us recall the occupancy scheme. Let $\mathcal{I}$ be a countable set and $\mathbf{p} = (p_i : i \in \mathcal{I})$ be a probability measure on $\mathcal{I}$. The occupancy scheme on $(\mathcal{I}, \mathbf{p})$ is described as follows. For all $i \in \mathcal{I}$ such that $p_i \neq 0$, one places a box at $i$. One then throws successively and independently $n$ balls in the boxes by assuming that each ball has probability $p_i$ of falling into the box located at $i$. We may be interested in the number of boxes containing exactly $j$ balls, or in the number of occupied boxes, etc.

We consider here a variant of the occupancy scheme which corresponds to a nested family of boxes. We introduce the infinite genealogical tree $\mathcal{U}$ :
$$
\mathcal{U} := \bigcup_{k = 0}^{\infty} \mathbb{N}^k,
$$
with $\mathbb{N} = \{1, 2, \dots \}$ and the convention $\mathbb{N}^0 = \{ \varnothing \}$. The elements of $\mathcal{U}$ are called individuals, and for every integer $k \in \mathbb{Z}_+$, the $k$-th generation of $\mathcal{U}$ is formed by the individuals $i \in \mathbb{N}^k$ (we write $|i| = k$). For every individual $i = (i_1, \dots, i_k)$ of $\mathcal{U}$ and $j \in \mathbb{N}$, the individual $i j = (i_1, \dots, i_k, j)$ is called the $j$-th child of $i$ and $i$ is the parent of $i j$. We suppose a real number $p_i \in [0,1]$ is assigned to each individual $i$, such that $p_{\varnothing} = 1$ and  $\sum_{j \in \mathbb{N}} p_{i j} = p_i$ for all $i \in \mathcal{U}$. Note that for all $k$, $\mathbf{p}(k) := (p_i : |i|=k)$ is a probability measure on $\mathbb{N}^k$. We can couple the occupancy schemes on $(\mathbb{N}^k, \mathbf{p}(k))$ as follows. Each individual $i$ should be viewed as a box, such that the box $i j$ is contained into the box $i$ for every $i \in \mathcal{U}$ and $j \in \mathbb{N}$ such that $p_{i j} \neq 0$. Initially the $n$ balls are thrown in the box located at $\varnothing$. Then one places successively and independently the $n$ balls in the boxes of the first generation by assuming that each ball has probability $p_i$ of falling into the box situated at $i \in \mathbb{N}$. Likewise, by iteration, a ball located in the box $i$ is placed independently of the others into the sub-box $i j$ with probability $p_{i j} / p_i$.

We denote by $H_{n,j}$ the first generation at which all the boxes have less than $j \geq 2$ balls when  $n$ have been thrown ($H_{n,j}$ is called a \emph{height}) and by  $G_{n, j}$ the first level at which there exists a box containing less than $j \geq 1$ balls ($G_{n,j}$ is called a \emph{saturation level}). Our aim is to study the asymptotic behaviours of $H_{n,j}$ and $G_{n,j}$ as $n$ tends to infinity when we consider a certain randomized version of $(p_i, i \in \mathcal{U})$.

More precisely, in the present work, we shall assume that we are given a random probability measure $\rho = (\rho_1, \rho_2, \dots)$ on $\mathbb{N}$. We assign to each individual $i$ an independent copy $\rho(i)$ of $\rho$. The real numbers $p_i$, $i \in \mathcal{U}$, are defined by induction~: $p_{\varnothing} := 1$ and $p_{i j} := p_i \rho_j(i)$, for all $i \in \mathcal{U}$ and $j \in \mathbb{N}$. It means that $\rho(i)$ describes how the mass of the individual $i$ is splitted to its chidren. This model is called the \emph{occupancy scheme of multiplicative cascades}. We put emphasis on the fact that there are two levels of randomness in our model, namely the arrangement of boxes  and the way one throws balls.

In the particular case when $\rho$ is supported by a finite number of integers, in the sense that $\# \{j : \rho_j > 0 \} \leq b$ a.s.\,\,for some integer $b \geq 2$, the height $H_{n,j}$ has a natural interpretation in terms of a special class of random split trees which have been considered \emph{e.g.}\,\,by Devoye \cite{devroko}. Specifically, imagine that each box has a rupture threshold of $j$, in the sense that when a ball falls into some box $i$ already containing $j-1$ balls, then this box is removed and the $j$ balls are shared out amongst the children of $i$ according to the random probability $\rho(i)$ (\emph{i.e.}\,\,conditionally on $\rho(i)$, each ball is put in the box $i j$ with probability $\rho_j(i)$, independently of the other balls). This procedure yields a random tree where all balls are stored at leaves, and $H_{n,j}$ is the height of this tree when $n$ balls have been thrown.

We further point out that the height $H_{n,j}$ also arises as a natural shattering time in homogeneous fragmentation chains, a class of partition-valued Markov chains. More precisely, the shattering time is defined as the first instant when all the blocks of the partition process have cardinality less than $j$. See \cite{bertoincoa} for background and Section~3.2 in \cite{bertocasc} for a description which is closely related to the present work. In a different direction, we mention the work of Haas \emph{et al.}\,\,\cite{hamipiwi} who associate another random tree to homogeneous fragmentation processes.

We shall show the following result~: there exist an integer $j^{\ast} \in \{ 2, 3, \dots \} \cup \{ \infty \}$ and a sequence of positive real numbers $C_2 >  \dots > C_{j^{\ast}} = C_{j^{\ast}+1} = \dots =: C^{\ast}$  such that for every integer $j \geq 2$, 
$$
H_{n,j} \underset{n \rightarrow \infty}{\sim} C_j \ln n \quad \textrm{a.s.}
$$
It should not be surprising that the heights $H_{n,j}$ have logarithmic asymptotics, as a large class of trees have such a behaviour, like random split trees.  On the other hand, it is remarkable that there exists a critical parameter $j^{\ast} \in \{ 2, 3, \dots \} \cup \{ \infty \}$ at which a phase transition occurs (provided that $2 < j^{\ast} < \infty$). It is a peculiar behaviour for a random split tree. Devroye indeed proved in \cite{devroko} that for a large class of random split trees such that their internal nodes contain at least one ball\,\footnote{We stress that this assumption is crucial in the proof of Theorem 1 in \cite{devroko} and seems to have been overlooked in the statement.}, then all the heights, regardless of the value of the rupture threshold $j$, have the same asymptotic behaviour (in probability). In particular, no phase transition occurs.

It may also be interesting to consider the situation where the parameter $j$ depends on the number of balls $n$. In the case of power functions, we shall prove that under some technical conditions that will be detailed below,
$$
H_{n, n^{\alpha}} \underset{n \rightarrow \infty}{\sim} (1-\alpha) C^{\ast} \ln n \quad \textrm{a.s.},
$$
for all $\alpha \in (0,1)$, so that there is no phase transition in the asymptotics of $H_{n, n^{\alpha}}$.

Concerning the asymptotics of the saturation level $G_{n,j}$, we shall prove the following result{~}: under some technical conditions that will be detailed below, there exists a constant $C_{\ast} < C^{\ast}$ such that for every positive integer $j$,
$$
G_{n,j}  \underset{n \rightarrow \infty}{\sim}  C_{\ast} \ln n  \quad \textrm{a.s.}
$$
In particular, there is no phase transition in the asymptotics of $G_{n,j}$. We shall also show that under some further technical conditions, for all $\alpha \in (0,1)$,
$$
G_{n, n^{\alpha}} \underset{n \rightarrow \infty}{\sim} (1-\alpha) C_{\ast} \ln n \quad \textrm{a.s.},
$$
so that there is also no phase transition in the asymptotic behaviours of $G_{n, n^{\alpha}}$.

Our approach essentially relies on the theory of branching random walks, and more precisely on their large deviations behaviours whose descriptions are due to Biggins \cite{biggins}. The construction of the boxes given by the multiplicative cascades indeed enables us to define a branching random walk giving the sizes of the boxes at each generation. We shall see that the critical parameter $j^{\ast}$ and the real numbers $C_2, \dots, C_{j^{\ast}}$ and $C_{\ast}$ are described by that branching random walk. Another key technique is Poissonization; instead of throwing exactly $n$ balls, one throws $\mathcal{P}_n$ balls, where $\mathcal{P}_n$ is a Poisson variable with parameter $n$ which is independent of $(\rho(i), i \in \mathcal{U})$. For every integer $k$, conditionally on the sizes of the boxes of the $k$-th generation, the numbers of balls per box of the $k$-th generation are thus independent Poisson variables.

The results will be stated in Section~\ref{formures}. We shall see the main techniques in Section~\ref{prelmin}. Section \ref{sectup}  will be devoted to the study of the heights. We shall first turn our attention to the upper bound. Due to the phase transition, we shall give two different proofs to show the lower bound of $H_{n,j}$, depending on whether the integer $j$ is less than the critical parameter $j^{\ast}$ or not. The results on the saturation levels will be proved in Section~\ref{satusect}.  Finally, in Section~\ref{studexap}, we shall explain heuristically why a phase transition may occur in the asymptotics of $H_{n,j}$, but not in those of $H_{n, n^{\alpha}}$, $G_{n,j}$ and $G_{n, n^{\alpha}}$.

\section{Formulation of the main results} \label{formures}

Recall that $\rho$ is a random probability measure on $\mathbb{N}$. We denote its law by $\nu$. If we denote by $\operatorname*{Prob}_{\mathbb{N}}$ the space of probability measures on $\mathbb{N}$, $\nu$ is a probability measure on $\operatorname*{Prob}_{\mathbb{N}}$, called the splitting law. We assume that $\nu$ is not geometric\,\footnote{Working with a geometric splitting law would induce a phenomenon of periodicity which we shall not discuss here for simplicity. However results similar to those proven in this work can be established by the same techniques for geometric splitting laws.}, in the sense that there is no real number $r \in (0,1)$ such that with probability one, all the masses of the atoms of $\rho$ belong to $\{ r^n, n \in \mathbb{Z}_+ \}$. In particular, the degenerate case when $\rho$ is a Dirac point mass a.s.\,\,is therefore excluded.

As explained in the introduction, we consider a family $(\rho(i), i \in \mathcal{U})$ of independent copies of $\rho$ labeled by the individuals of the genealogical tree $\mathcal{U}$. The multiplicative cascade construction defines for each generation $k$ a probability measure $(p_i, |i| = k)$ on $\mathbb{N}^k$. Taking logarithm of masses, we may encode the latter by the following random point measure on $\mathbb{R}_+$
$$
Z^{(k)} (\textrm{d} y) := \sum_{i : |i| = k} \delta_{- \ln p_i} (\textrm{d} y),
$$
where $\delta_z$ stands for the Dirac point mass at $z$. Note that if $p_i = 0$, \emph{i.e.}\,\,if there is no box at $i$, the individual $i$ is omitted in the sum defining $Z^{(k)}$. Likewise, in the sequel, we shall always consider individuals which have a positive mass.  It follows immediately from the structure of the multiplicative cascades that $(Z^{(k)}, k \in \mathbb{Z}_+)$ is a branching random walk, in the sense that for every integers $k, k' \geq 0$, $Z^{(k+k')}$ is obtained from $Z^{(k)}$ by replacing each atom $z$ of $Z^{(k)}$ by a family $\{ z+y, y \in \mathcal{Y}_z \}$, where each $\mathcal{Y}_z$ is an independent copy of the family of atoms of $Z^{(k')}$.

Let us introduce quantities defined via the splitting law $\nu$. First, we define the Laplace transform of the intensity measure $Z^{(1)}$ by
$$
\textrm{L}(\theta) := \mathbb{E} \left [ \left \langle Z^{(1)}, e^{-\theta \cdot} \right \rangle \right ]
$$
for $\theta \in \mathbb{R}$. We can also write
$$
\textrm{L}(\theta) = \mathbb{E} \left [  \sum_{j \in \mathbb{N}} \rho_j^{\theta} \right ] = \int_{\operatorname*{Prob}_{\mathbb{N}}} \left ( \sum_{j \in \mathbb{N}} p_j^{\theta} \right ) \nu(\textrm{d} \mathbf{p})
$$
with the convention that $p^{\theta} = 0$ when $p=0$ even when $\theta \leq 0$. Because  $\rho$ is not a Dirac point mass a.s.,\,\,$\textrm{L}(0) > 1$. The function $\textrm{L} : \mathbb{R} \rightarrow (0, \infty]$ is decreasing with $\textrm{L}(1) = 1$. We define
$$
\underline{\theta} := \inf \{  \theta \in \mathbb{R} : \textrm{L}(\theta) < \infty \},
$$
so that $\textrm{L}(\theta) < \infty$ when $\theta > \underline{\theta}$. Note that $\underline{\theta}$ may be negative. We can show from Hölder's inequality that $\ln \textrm{L}$ is a convex function, which implies the convexity of $\operatorname*{L}$ and that the function
$$
\varphi : \theta \mapsto \ln \textrm{L} (\theta) - \theta \frac{\textrm{L}' (\theta)}{\textrm{L} (\theta)}
$$
is increasing on $(\underline{\theta}, 0)$ and decreasing on $(0, \infty)$. As $\textrm{L}(1)=1$ and $\textrm{L}$ decreases, we have $\varphi(1) = - \textrm{L}'(1) > 0$, and thus the set of $\theta \in (\underline{\theta}, \infty)$ such that $\varphi(\theta) > 0$ is a non-empty open interval $(\theta_{\ast}, \theta^{\ast})$~:
$$
\theta_{\ast} := \inf \{ \theta > \underline{\theta} : \varphi(\theta) > 0   \}  \quad \textrm{and} \quad \theta^{\ast} := \sup \{ \theta > \underline{\theta} : \varphi(\theta) > 0   \}.
$$
Note that $\theta^{\ast} > 1$ and $\theta_{\ast} < 0$ if $\underline{\theta} < 0$. The convexity of $\ln \rm{L}$ implies that $-\textrm{L}/\textrm{L}'$ is an increasing function. We define
$$
C_{\ast} := \lim_{\theta \downarrow \theta_{\ast}} \downarrow - \frac{\textrm{L}(\theta)}{\textrm{L}'(\theta)} \quad \textrm{and} \quad C^{\ast} := \lim_{\theta \uparrow \theta^{\ast}} \uparrow - \frac{\textrm{L}(\theta)}{\textrm{L}'(\theta)}.
$$

Inspired by the article of Hu and Shi in \cite{hushi} (see Lemma \ref{thelimsuphushi} below), we shall sometimes need the following assumption : there exists $\delta > 0$ such that
 \begin{equation} \label{deuxhypohushi}
\operatorname*{L}(- \delta) < \infty \quad \textrm{and} \quad  \int_{\operatorname*{Prob}_{\mathbb{N}}} \left ( \sum_{i \in \mathbb{N}} \textbf{1}_{p_i > 0}  \right )^{1 + \delta} \nu (\textrm{d} \mathbf{p}) < \infty.
 \end{equation}
To study the asymptotics of the saturation levels, we shall sometimes need the following hypothesis{~}:
\begin{equation} \label{troishypohushi}
- \infty < \theta_{\ast} < 0 \quad \textrm{and} \quad \varphi(\theta_{\ast}) = 0.
\end{equation}

We can now state the results that we shall prove. Concerning the asymptotic behaviours of the heights $H_{n,j}$, we have the following results which complete and improve Proposition~2 in \cite{bertocasc}.
\begin{theoreme} \label{coucoudidi}
Let $j \geq 2$ be an integer, set
$$
C_j := \left \lbrace
\begin{array}{ll}
- j/  \ln \operatorname*{L}(j) & \textrm{if } j < \theta^{\ast}, \\
C^{\ast}  & \textrm{if } j \geq \theta^{\ast}.
\end{array}
\right. 
$$
Then
$$
H_{n,j} \sim C_j \ln n \quad \textrm{a.s.}
$$
More precisely,
\begin{equation} \label{theoprenon}
H_{n,j} \leq C_j \ln n + O (\ln \ln  n) \quad \textrm{a.s.},
\end{equation}
and (\ref{theoprenon}) is an equality if $j < \theta^{\ast}$ or if (\ref{deuxhypohushi}) holds.
\end{theoreme}

We see that there is a phase transition in the asymptotic behaviour of $H_{n,j}$ at the integer $\lceil \theta^{\ast} \rceil$ when $2 < \theta^{\ast} < \infty$. We point out that it may happen that $\theta^{\ast} \leq 2$ or $\theta^{\ast} = \infty$, in which case there is no phase transition. For instance, one can show that if $\rho = (\rho_1, \rho_2, \dots)$ with $\rho_1 = 1-0.75 U$, $\rho_j = 0.05 U$ for all $j \in {2, \dots, 16}$ and $\rho_j = 0$ for all $j \geq 17$, where $U$ is uniformly distributed on $(0,1)$, then $\theta^{\ast} < 1.99$. On the other hand, for every $\alpha \in [1/2, 1)$, if $\rho$ is equal to $(1/2, 1/2, 0, \dots)$ with probability $\alpha$ and to $(1/3, 1/3, 1/3, 0, \dots)$ with probability $1-\alpha$, then $\theta^{\ast} = \infty$.

In this direction, we also point out that the critical parameter $\theta^{\ast}$ is finite whenever
$$
\| \max_{j \in \mathbb{N}} \rho_j \|_{\infty} = 1,
$$
where $\rho = (\rho_j, j \in \mathbb{N})$ denotes a random probability measure on $\mathbb{N}$ with law $\nu$. Indeed, it is easily seen that 
$$
\lim_{\theta \rightarrow \infty} \operatorname*{L}(\theta)^{1 / \theta} = \| \max_{j \in \mathbb{N}} \rho_j \|_{\infty}.
$$
If the right-hand side equals 1, then $g : \theta \mapsto - \ln \operatorname*{L}(\theta) / \theta$ has limit 0 at infinity. Now $g$ is derivable and $g(1) = 0$, so there exists $\theta_{0} \in (1, \infty)$ such that $g'(\theta_0) = 0$. As $g'(\theta) = \theta^{-2} \varphi(\theta)$, we conclude that $\theta_0 = \theta^{\ast} < \infty$.

We may also be interested in the asymptotics of $H_{n,n^{\alpha}}$, where $\alpha \in (0,1)$. We shall prove the following result.
\begin{prop} \label{mscjetaime}
Suppose $\theta^{\ast} < \infty$.  Let $\alpha \in (0,1)$. 
Then
$$
H_{n,n^\alpha} \sim (1- \alpha) C^{\ast} \ln n \quad \textrm{a.s.}
$$
Furthermore,
\begin{equation} \label{theodeenon}
H_{n, n ^{\alpha}} \leq (1- \alpha) C^{\ast} \ln n + O (\ln \ln n) \quad \textrm{a.s.},
\end{equation}
and (\ref{theodeenon}) is an equality whenever (\ref{deuxhypohushi}) holds.
\end{prop}

\begin{remark} \rm
The final assertions in Theorem~\ref{coucoudidi} and Proposition~\ref{mscjetaime} rely on the work of Hu and Shi \cite{hushi}. McDiarmid's setting in \cite{mcdiarmid} can however be considered; we can prove that the results stated in Theorem~\ref{coucoudidi} and in Proposition~\ref{mscjetaime} still hold if the assumption (\ref{deuxhypohushi}) is replaced by the following~:
 \begin{equation} \label{hypomcdiarmid}
\int_{\operatorname*{Prob}_{\mathbb{N}}} \left ( \sum_{i \in \mathbb{N}} \textbf{1}_{p_i > 0}  \right )^{2} \nu (\textrm{d} \mathbf{p}) < \infty.
 \end{equation}
 
For instance, suppose that $\rho = (\rho_1, 1-\rho_1, 0, 0, \dots)$, where $\rho_1$ is a random variable with density $\textbf{1}_{0 < x < \operatorname*{e}^{-1}} x^{-1} \ln^{-2} x \textrm{d} x$. Then $\underline{\theta} = 0$, so  (\ref{deuxhypohushi}) does not hold. Nonetheless, $\rho$ is supported by two integers a.s., so (\ref{hypomcdiarmid}) holds. As a result, for all $j \geq 2$, $H_{n,j} = C_j \ln n +O (\ln \ln n)$ a.s. Furthermore, as $\| \max_{j \in \mathbb{N}} \rho_j \|_{\infty} = 1$, the discussion below Theorem~\ref{coucoudidi} ensures that $\theta^{\ast}$ is finite, so for all $\alpha \in (0,1)$, $H_{n,n^{\alpha}} = (1- \alpha) C^{\ast} \ln n +O (\ln \ln n)$ a.s.

For the sake of simplicity, we shall show Theorem~\ref{coucoudidi} and Proposition~\ref{mscjetaime} only in the setting of \cite{hushi} (see Proposition~\ref{lowerno1poup} below). The general proof can however be easily carried out.
\end{remark}

Concerning the asymptotic behaviours of the saturation levels $G_{n,j}$, we shall prove the following theorem.
\begin{theoreme} \label{maisonneuve}
Let $j \geq 1$ be an integer.
\begin{itemize}
\item If $\theta_{\ast} = -\infty$, then
$$
G_{n,j} \sim C_{\ast} \ln n \quad \textrm{a.s.}
$$
\item Suppose that (\ref{troishypohushi}) holds. Then 
$$
G_{n,j} \sim C_{\ast} \ln n \quad \textrm{a.s.}
$$
More precisely,
\begin{equation} \label{theoprenonkiooo}
G_{n,j} \geq C_{\ast} \ln n + O (\ln \ln  n) \quad \textrm{a.s.},
\end{equation}
and (\ref{theoprenonkiooo}) is an equality if (\ref{deuxhypohushi}) holds.
\end{itemize}
\end{theoreme}

A sufficient condition to guarantee that  $- \infty < \theta_{\ast} < 0$ and $\varphi ( \theta_{\ast} ) = 0$ is :
$$
\underline{\theta} > - \infty \quad \textrm{and} \quad \lim_{\theta \downarrow \underline{\theta}} \uparrow \operatorname*{L} (\theta) = \infty.
$$
Indeed, suppose this condition fulfilled. Imagine that for all $\theta < 0$, $\varphi (\theta) > 0$. Let $\theta_0 \in (\underline{\theta}, 0)$. As $\theta \mapsto \ln \operatorname*{L}(\theta) / \theta$ is decreasing (its derivative is $-\theta^{-2} \varphi(\theta) < 0$), we have for all $\theta < \theta_0${~}:
$$
\ln \operatorname*{L} (\theta) \leq \theta \frac{\ln \operatorname*{L} (\theta_0)}{\theta_0} \leq \underline{\theta} \frac{\ln \operatorname*{L} (\theta_0)}{\theta_0},
$$
which contradicts $\lim_{\theta \rightarrow \underline{\theta}} \operatorname*{L} (\theta) = \infty$.

We shall also show the following proposition.
\begin{prop} \label{laroseara}
Suppose that (\ref{troishypohushi}) holds. Let $\alpha \in (0,1)$. Then
$$
G_{n,n^{\alpha}} \sim (1- \alpha) C_{\ast} \ln n \quad \textrm{a.s.}
$$
Moreover,
\begin{equation} \label{theoprenonkioooij}
G_{n,n^{\alpha}} \geq (1-\alpha) C_{\ast} \ln n + O (\ln \ln  n) \quad \textrm{a.s.},
\end{equation}
and (\ref{theoprenonkioooij}) is an equality if (\ref{deuxhypohushi}) holds.
\end{prop}

We now conclude this section by discussing an illustrative example. Consider the case $\rho = (U, 1-U, 0, 0, \dots)$, where $U$ is uniformly distributed on $[0,1]$. Then $\operatorname*{L}(\theta) = 2/(\theta +1)$, so $\underline{\theta} = -1$ and (\ref{deuxhypohushi}) holds.  The discussions below the theorems yield~: $\theta^{\ast} < \infty$ and (\ref{troishypohushi}) hold. Hence, all our results may be applied. Easy calculations yield $\varphi(\theta) = \ln 2 - \ln (\theta +1) + \theta/(\theta+1)$ and $\lceil \theta^{\ast} \rceil = 4$, so that there is a phase transition. We can show that $C_2 = 2 / \ln (3/2) \approx 4,93260...$, $C_3 = 3 /\ln 2 \approx 4,32808...$ and that $C^{\ast} > C_{\ast}$ are the solutions of the equation 
$$
\ln \left (  \frac{2}{c} \right ) + \frac{c-1}{c} = 0,
$$
\emph{i.e.}\,\,$C^{\ast} \approx 4,31107...$ and $C_{\ast} \approx 0,37336...$ 

As $\rho$ is supported by two integers, our model may be interpreted in terms of random split trees; the procedure described in the introduction yields a random tree where all balls are stored at leaves. We denote by $T_{n,j}$ the tree obtained when $n$ balls have been thrown and when the boxes have a rupture threshold of $j$. We now define another random split tree $\tilde{T}_{n,j}$ also related to our model. We imagine that each box has a rupture threshold of $j$, but when a ball falls into some box $i$ already containing $j-1$ balls, it remains in that box and the $j-1$ other balls are shared out amongst the two children of $i$ independently and with probability $(1/2, 1/2)$. No other ball is then allowed to be stored at the box $i$~: when a ball falls into $i$, it is placed into one of the two children of $i$ with probability $(1/2, 1/2)$, and one has to consider whether the ball stays at that box or not. We denote by $\tilde{T}_{n,j}$ the tree obtained when $n$ balls have been thrown. Note that each internal node of $\tilde{T}_{n,j}$ has exactly one ball, whereas its leaves have less than $j$ balls. We see that the height $\tilde{H}_{n,j}$ of the tree $\tilde{T}_{n,j}$ is less than or equal to the height $H_{n,j}$ of the tree $T_{n,j}$.

It is remarkable that $\tilde{T}_{n,2}$ has the law of the random binary search tree. Robson \cite{robson} was the first to be interested in the height of the random binary search tree. In \cite{pittel}, Pittel studied its height and its saturation level. Devroye proved in \cite{devro86, devro87} that the saturation level is asymptotically equivalent to $C_{\ast} \ln n$ in probability as in our model, but that the height $\tilde{H}_{n,2}$ is equivalent to $C^{\ast} \ln n$ in probability, so that $\tilde{H}_{n,2}$ is not equivalent to $H_{n,2}$. Regarding the random binary search tree as a random split tree, he proved in \cite{devroko} that all the heights $\tilde{H}_{n,j}$ (regardless of the value of $j$) of the random binary search tree are asymptotically equivalent to $C^{\ast} \ln n$ in probability. In particular, there is no phase transition, contrary to our case. There is of course no contradiction, as the random binary search tree does not correspond to a model treated in this work. Indeed, in terms of random split trees, internal nodes of the binary search tree retain exactly one ball whereas all the balls are stored at leaves in our case.

\section{Preliminaries} \label{prelmin}

\subsection{Some results on branching random walks}

In this section, we recall some results on branching random walks that will be used in the proofs.

We begin by stating a key result obtained by Biggins in \cite{biggins}. For every $\theta > \underline{\theta}$, we introduce
$$
W^{(k)} (\theta) := \operatorname*{L} (\theta)^{-k} \left \langle  Z^{(k)}, e^{- \theta \cdot }       \right \rangle =  \operatorname*{L}(\theta)^{-k} \sum_{i : |i| = k} p_i^{\theta}.
$$
For every $k \in \mathbb{Z}_+ \cup \{ \infty \}$, we denote by $\mathcal{F}_k$ the $\sigma$-algebra generated by $(p_i, |i| \leq k)$.
\begin{lemme} \label{martibigg}
For every $\theta > \underline{\theta}$, $(W^{(k)}(\theta), k \in \mathbb{Z}_+)$ is a martingale with respect to the filtration $\left (\mathcal{F}_k \right )_{k \in \mathbb{Z}_+}$. Moreover, if $\theta \in (\theta_{\ast}, \theta^{\ast})$, it is bounded in $L^{\gamma}(\mathbb{P})$ for some $\gamma > 1$ and therefore uniformly integrable, and its terminal value 
$$
W(\theta) := \lim_{k \rightarrow \infty} W^{(k)}(\theta)
$$
is positive a.s.
\end{lemme}

Applying Corollary 4 in \cite{biggins}, we can obtain a precise estimate of the number of boxes at generation $k$ with size of order $\exp (k \operatorname*{L}'(\theta) / \operatorname*{L}(\theta))$. Specifically
\begin{lemme} \label{martibigg2}
For all real numbers $a$ and $b$ such that $a<b$ and for all $\theta \in (\theta_{\ast}, \theta^{\ast})$, we have with probability one that
$$
\lim_{k \rightarrow \infty}
\sqrt{k} e^{-k \varphi(\theta)} \# \left \{ i \in \mathbb{N}^k :  \exp \left ( -b + k {\operatorname*{L}}' (\theta) / \operatorname*{L}(\theta) \right )  \leq p_i \leq  \exp \left ( -a + k {\operatorname*{L}}'(\theta) / \operatorname*{L}(\theta) \right )  \right \}
$$
$$
= \frac{1}{\sqrt{2 \pi}} \frac{e^{\theta b} - e^{\theta a}}{\theta} \left ( \frac{\operatorname*{L}''(\theta)}{\operatorname*{L}(\theta)}  - \left (  \frac{\operatorname*{L}'(\theta)}{\operatorname*{L}(\theta)} \right )^2  \right )^{-1/2} W(\theta) .
$$
\end{lemme}

We next turn our interest to the asymptotic behaviours of extreme sizes of boxes at generation $k${~}:
$$
\underline{p}(k) := \inf \left \{ p_i : p_i > 0, |i| = k \right \} \quad \textrm{and} \quad \overline{p}(k) := \sup \left \{ p_i :  |i| = k \right \} .
$$
We have the general following results :
\begin{lemme} \label{martibigg3}
If $\theta^{\ast} < \infty$, then
$$
\lim_{k \rightarrow \infty} \left ( - \ln \overline{p}(k) - k / C^{\ast} \right ) = \infty \quad \textrm{and} \quad \lim_{k \rightarrow \infty} - \frac{\ln \overline{p}(k)}{k} = \frac{1}{C^{\ast}} \quad \textrm{a.s.}
$$
Likewise, if (\ref{troishypohushi}) holds, then
$$
\lim_{k \rightarrow \infty}  \left ( \ln \underline{p}(k) + k / C_{\ast} \right ) = \infty \quad \textrm{and} \quad \lim_{k \rightarrow \infty} - \frac{\ln \underline{p}(k)}{k} = \frac{1}{C_{\ast}} \quad \textrm{a.s.}
$$
\end{lemme}
To have precise estimates of the behaviours of the heights and of the saturation levels, we shall sometimes need sharper results of the asymptotics of $\overline{p}(k)$ and $\underline{p}(k)$. The following result, proved by Hu and Shi in \cite{hushi}, will be very useful.
\begin{lemme} \label{thelimsuphushi} We assume that (\ref{deuxhypohushi}) holds.
\begin{itemize}
\item If $\theta^{\ast} < \infty$, then 
$$
\limsup_{k \rightarrow \infty} \frac{- \ln \overline{p}(k) - k / C^{\ast}}{\ln k} =  \frac{3}{2 \theta^{\ast}} \quad \textrm{a.s.}
$$
\item If (\ref{troishypohushi}) holds, then
$$
\limsup_{k \rightarrow \infty}  \frac{\ln \underline{p}(k) + k / C_{\ast}}{\ln k} = -\frac{3}{2 \theta_{\ast}} \quad \textrm{a.s.}
 $$
\end{itemize}
\end{lemme}
Addario-Berry and Reed in \cite{addarioreed} also studied the minima in branching random walks, but they require a stronger condition than (\ref{deuxhypohushi}). Their assumption is however fulfilled for random split trees.

\subsection{Poissonization}

Let us present the methods used in the proofs. We have to consider the number of balls belonging to each box at each generation. Now, conditionally on $\mathcal{F}_k$, when $n$ balls have been thrown, the number of balls in the box $i$, $|i| = k$, follows a binomial law of parameter $n p_i$. Furthermore, these random variables are not independent. A classical idea to circumvent those difficulties (see for instance Gnedin \emph{et al.}\,\,in \cite{gnhapi} or Holst in \cite{holst}) is to consider a randomized version of the total number of balls{~}: instead of throwing initially $n$ balls, one throws $\mathcal{P}_n$ balls, where $\mathcal{P}_n$ is a Poisson variable with parameter $n$ which is independent of $(\rho(i), i \in \mathcal{U})$.

More precisely, we suppose we are given a standard Poisson process $\left ( \mathcal{P}_x \right )_{x \geq 0}$  independent of $\mathcal{F}_{\infty}$. For every individual $i \in \mathcal{U}$ and for every $x \in (0, \infty)$, we denote by $\mathcal{C}(i;x)$ the number of balls at $i$ when the first $\mathcal{P}_x$ balls have been thrown. For all $x, y \in (0, \infty)$ and $k \in \mathbb{Z}_+$, we denote by $\mathcal{N}_{x, y}(k)$ the number of boxes at generation $k$ containing at least $y$ balls when the first $\mathcal{P}_x$ balls have been thrown :
$$
\mathcal{N}_{x,y}(k) := \# \{ i \in \mathbb{N}^k : \mathcal{C}(i ; x) \geq y    \}.
$$ 
Conditionally on $\mathcal{F}_k$, the random variables $(\mathcal{C}(i ; n))_{|i|=k}$ are independent Poisson variables with parameters $n p_i$. That is why in our proofs, we shall first focus on $\mathcal{N}_{n,j}(k)$. We shall then show that $\mathcal{N}_{n,j}(k)$ is close to the number $N_{n,j}(k)$ of boxes at generation $k$ containing at least $j$ balls when exactly $n$ balls have been thrown.
Similarly, we denote by $\mathcal{M}_{x, y}(k)$ the number of boxes at generation $k$ containing less than $y$ balls when the first $\mathcal{P}_x$ balls have been thrown :
$$
\mathcal{M}_{x,y}(k) := \# \{ i \in \mathbb{N}^k : \mathcal{C}(i ; x) < y    \}
$$ 
and by $M_{n,j}$ the number of boxes at generation $k$ containing less than $j$ balls when $n$ balls have been thrown.

We now prove two estimates using the technique of Poissonization that will then be used in the sequel.
\begin{lemme} \label{lemmemodif}
Let $p \in (0, \infty)$. There exist two finite constants $c(p)$ and $d(p)$ such that
\begin{equation} \label{numbonelemmod}
\sup_{\scriptstyle j \geq p \atop \scriptstyle j \in \mathbb{N}}  \sup_{x > 0} \mathbb{P} \left (  \mathcal{P}_x \geq j  \right ) j^p x^{-p}  \leq c(p),
\end{equation}
and
\begin{equation} \label{numbtwolemmod}
\sup_{j \in \mathbb{N}}  \sup_{x > 0} \mathbb{P} \left (  \mathcal{P}_x < j  \right )  j^{-p}  x^{p} \leq d(p).
\end{equation}
\end{lemme}
\textbf{Proof :} We begin by showing (\ref{numbonelemmod}).  Let $j \geq p$ be an integer. Let $x>0$. By Markov's inequality, 
$$
\mathbb{P} \left (  \mathcal{P}_x \geq j  \right ) = \mathbb{P} \left ( \mathcal{P}_x^p \geq j^p \right ) \leq j^{-p} \mathbb{E} \left[   \mathcal{P}_x^p  \textbf{1}_{  \mathcal{P}_x \geq j} \right ]  .
$$
Therefore we only have to bound from above
$$
x^{-p} \mathbb{E} \left[ \mathcal{P}_x^p  \textbf{1}_{  \mathcal{P}_x \geq j  } \right] =   x^{-p} \sum_{k=j}^{\infty} k^p e^{-x} \frac{x^k}{k!}. 
$$ 
As $\Gamma(k-p+1) k^p / k! \rightarrow 1$ as $k$ tends to infinity, there exists a finite constant $c(p) \geq 1$ such that for all $k > p-1$, 
$$
\frac{k^p}{k!} \leq \frac{c(p)}{\Gamma(k-p+1)}.$$ As $j \geq p > p-1$, we thus have
$$
 x^{-p} \sum_{k=j}^{\infty} k^p e^{-x} \frac{x^k}{k!} \leq c(p) e^{-x} \sum_{k=j}^{\infty} \frac{x^{k-p}}{\Gamma(k-p+1)}  = c(p) e^{-x} \sum_{k=0}^{\infty} \frac{x^{k+u}}{\Gamma(k+u+1)},
$$
where $u:=j-p \geq 0$. Applying the formulae 6.5.1, 6.5.4 and 6.5.29 in \cite{handbook}, we get that 
$$
\sum_{k=0}^{\infty} \frac{x^{k+u}}{\Gamma(k+u+1)} \leq e^x,
$$
which proves (\ref{numbonelemmod}). 

To show (\ref{numbtwolemmod}), we write by Markov's inequality :
$$
\mathbb{P} (\mathcal{P}_x < j)  = \mathbb{P} \left ( \left ( \mathcal{P}_x + 1 \right )^{-p} \geq j^{-p} \right ) \leq j^{p} \mathbb{E} \left[ \left (  \mathcal{P}_x+1 \right )^{-p} \right ]  = j^{p} \sum_{k=0}^{\infty} (k+1)^{-p} e^{-x} \frac{x^k}{k!}.
$$
As $\Gamma(k+p+1) (k+1)^{-p} / k!  \rightarrow 1$ as $k$ tends to infinity, there exists a finite constant $d(p) \geq 1$ such that for all $k \in \mathbb{Z}_+$, $\Gamma(k+p+1) (k+1)^{-p} / k! \leq d(p)$. We thus have
$$
 \mathbb{P} (\mathcal{P}_x < j)  \leq d(p) j^p x^{-p} e^{-x} \sum_{k=0}^{\infty} \frac{x^{k+p}}{\Gamma(k+p+1)}.
$$
We conclude as before.
\hfill $\square$

\begin{remark}  \label{remstre} \rm
We stress that in (\ref{numbonelemmod}), the integer $j$ cannot be less than $p$. This restriction lies in the heart of the existence of the phase transition in the asymptotics of $H_{n,j}$ (see Proposition~\ref{deuxxpromsc} below, and more precisely Lemma \ref{lempoissomsc}). On the other hand, there is no restriction in (\ref{numbtwolemmod}); no phase transition appears in the asymptotics of $G_{n,j}$.
\end{remark}

\section{Study of the heights} \label{sectup}

In this section, we prove Theorem~\ref{coucoudidi} and Proposition~\ref{mscjetaime}. We are first interested in the upper bound. We shall then focus on the lower bound. We stress that the phase transition is glimpsed at the beginning of the study of the upper bound (see Remark \ref{genevibra} below). It is however proved in the paragraph dealing with the lower bound.

\subsection{Upper bound}

Equation (\ref{numbonelemmod}) will enable us to have a uniform upper bound of $N_{n, j}(k)$ independent of the sizes of the boxes that will eventually lead to the inequalities (\ref{theoprenon}) and (\ref{theodeenon}).

\begin{prop} \label{deuxxpromsc} Let $j \geq 1$ be an integer and $\alpha \in [0, 1)$ such that $(j, \alpha) \neq (1, 0)$. Then
$$
\limsup_{n \rightarrow \infty} \frac{H_{n,j n^ \alpha} - (1 - \alpha) \frac{\theta}{- \ln \operatorname*{L}(\theta)}  \ln n }{ \ln \ln n } \leq   \frac{1}{- \ln \operatorname*{L}(\theta)} \quad \textrm{a.s.}
$$
for every $\theta > 1$ if $\alpha > 0$, and for every $\theta \in (1, j]$ if $\alpha = 0$. 
\end{prop}

\begin{remark} \label{genevibra} \rm
 Suppose that $\alpha = 0$. We deduce from Proposition~\ref{deuxxpromsc} that{~}:
$$
H_{n,j} \leq \min_{\theta \in (1.j]} h(\theta) \ln n + O(\ln \ln n) \quad \textrm{a.s.},
$$
where $h$ is the function $\theta \mapsto - \theta / \ln \operatorname*{L}(\theta)$. Now, $h$ is derivable and $h'(\theta) = - \varphi(\theta) \ln^{-2} \operatorname*{L}(\theta)$, so $h$ is decreasing on $(1, \theta^{\ast}]$ and increasing on $[\theta^{\ast}, \infty)$. The minimum of $h$ is therefore
\begin{itemize}
\item $-j / \ln \operatorname*{L}(j)$ if $j < \theta^{\ast}$, 
\item $- \theta^{\ast} / \ln \operatorname*{L}(\theta^{\ast}) = C^{\ast}$ if $j \geq \theta^{\ast}$.
\end{itemize}
The phase transition is shown. It will be proved in the study of the lower bound.
\end{remark}

Recall that $H_{n,j n^{\alpha}} \leq k$ if and only if at genetation $k$, every box contains less than $j n^{\alpha}$ balls when $n$ balls have been thrown. Because we want to show that $H_{n, j n^{\alpha}}$ is bounded from above by  $\left ( (1 - \alpha) \frac{ \theta}{- \ln \operatorname*{L}(\theta)} + \varepsilon \right ) \ln(n)$, we take $k \approx \left (  (1 - \alpha) \frac{\theta}{-\ln \operatorname*{L}(\theta)} + \varepsilon \right ) \ln(n)$. In other words, one initially throws $n \approx \exp \left (  k \left (  \frac{1}{1 - \alpha} \frac{- \ln \operatorname*{L}(\theta)}{\theta}- \varepsilon \right ) \right )$ balls and we show that every box of the $k$-th generation contains less than $j n^{\alpha}$ balls.

Let $u$ and $u'$ be two real numbers such that
$$
u > \frac{1}{(1 - \alpha) \theta } \quad \textrm{and} \quad u < u' < \frac{1}{\alpha} \left ( u - \frac{1}{\theta}  \right ),
$$
where for $\alpha = 0$, the second condition reduces to $u' > u$.  Define for every $k \geq 1$
$$
x_k := k^{-u} \exp \left (  k \frac{1}{1-\alpha} \frac{- \ln \operatorname*{L}(\theta)}{\theta} \right ) \quad \textrm{and} \quad
\phi_k := k^{-u'} \exp \left ( k \frac{1}{1-\alpha} \frac{-\ln \operatorname*{L}(\theta)}{\theta} \right ).
$$ 
Informally, $x_k$ corresponds to the number of balls thrown when we consider boxes at generation $k$. Note that $x_k \approx \exp \left (  k \left (  \frac{1}{1 - \alpha} \frac{- \ln \operatorname*{L}(\theta)}{\theta}- \varepsilon \right ) \right )$.

As mentionned in the preliminaries, the argument relies on Poissonization. We are first interested in $\mathcal{N}_{x,y} (k)$ and we shall see how to depoissonize.

\begin{lemme} \label{lempoissomsc}
For almost all $\omega$, there exists $k_0(\omega)$ such that
$$
\mathcal{N}_{x_k,j \phi_k^{\alpha}} (k)= 0, \quad \textrm{for all $k \geq k_0(\omega)$}.
$$
\end{lemme}
\textbf{Proof :} Let $x > 0$ and $y \geq  \theta$. We calculate $\mathbb{E}\left[ \mathcal{N}_{x,y} (k)|\mathcal{F}_k \right] $. We write
$$
\mathcal{N}_{x,y} (k) = \sum_{i : |i|=k} \mathbf{1}_{\mathcal{C}(i;x) \geq y} = \sum_{i : |i|=k} \mathbf{1}_{\mathcal{C}(i;x) \geq \lceil y \rceil}.
$$
Conditionally on $\mathcal{F}_k$, $(\mathcal{C}(i;x))_{|i|=k}$ are Poisson variables with parameters $x p_i$, so
$$
\mathbb{E}\left[ \mathcal{N}_{x,y}(k)|\mathcal{F}_k \right]  =  \sum_{i : |i|=k} \mathbb{P}(\mathcal{P}_{x p_i} \geq  \lceil y \rceil).
$$
As $\lceil y \rceil \geq \theta$, (\ref{numbonelemmod}) ensures that
$$
\mathbb{E}\left[ \mathcal{N}_{x,y}(k)|\mathcal{F}_k \right]  \leq  \sum_{i : |i|=k} c( \theta) \lceil y \rceil ^{-\theta} (x p_i)^{\theta}
$$
(see Remark \ref{remstre}). Now, as $y \geq \theta$,
$$
 \lceil y \rceil    \geq  \frac{ y - 1 }{y}  y \geq \frac{ \theta - 1}{\theta}  y.
$$
In the notations of Lemma \ref{martibigg}, we have :
$$
\mathbb{E}\left[ \mathcal{N}_{x,y}(k)|\mathcal{F}_k \right]  \leq c'(\theta)  y^{-\theta} x^{\theta} \operatorname*{L}(\theta)^{k} W^{(k)} ( \theta),
$$
where $c'(\theta) := c(\theta) (1 - 1 / \theta)^{- \theta}$. We finally get 
$$
\mathbb{E}\left[ \mathcal{N}_{x,y}(k) \right ]  \leq c'(\theta)  y^{-\theta} x^{\theta} \operatorname*{L}(\theta)^{k}.
$$
Taking $ x = x_k$ and $y = j \phi_k^{\alpha}$, we get for all $k$ sufficiently large so that $j \phi_k^{\alpha} \geq \theta$ (note that if $\alpha = 0$, then for all $k$, $j \phi_k^{\alpha} \geq \theta$) :
$$
\mathbb{E} \left[ \mathcal{N}_{x_k,j \phi_k^{\alpha}}(k) \right]   \leq  c'(\theta) j^{-\theta}k^{\theta(u' \alpha - u)}.
$$
Now $\theta ( u' \alpha - u) < -1$. As a consequence
$$
\mathbb{E} \left[  \sum_{k \in \mathbb{N}}  \mathbf{1}_{\mathcal{N}_{x_k,j \phi_k^{\alpha}}(k)  \geq 1}  \right] \leq  \mathbb{E} \left[  \sum_{k \in \mathbb{N}}  \mathcal{N}_{x_k,j \phi_k^{\alpha}}(k) \right] < \infty. 
$$
In particular, there is an a.s.\,\,finite number of integers $k$ such that $\mathcal{N}_{x_k,j \phi_k^{\alpha}}(k)  \geq 1$.
\hfill $\square$

\ \\
We now show how to have information on $N_{n, j n^{\alpha}} (k)$ itself. 

\ \\
\textbf{Proof of Proposition~\ref{deuxxpromsc} :} Let $E_k$ denote the event $ \left \{  \mathcal{P}_{x_k} \leq \phi_{k+1}  \right \}$.  Because $u' > u$, $\phi_{k+1} / x_k$ tends to 0 as $k$ tends to infinity, so there exists an integer $k_1$ such that for all $k \geq k_1$, $\phi_{k+1} \leq x_k / 2$. Consequently, 
$$
\mathbb{P} \left( E_k\right)  \leq   \mathbb{P} \left(  \mathcal{P}_{x_k} \leq x_k/2 \right) = \mathbb{P} \left(  \mathcal{P}_{x_k} - x_k \leq - x_k/2 \right) \leq  \mathbb{P} \left( | \mathcal{P}_{x_k} - x_k| \geq  x_k/2 \right).
$$
The variance of the Poisson variable $\mathcal{P}_{x_k}$ being $x_k$, we get by Chebichev's inequality:
$
\mathbb{P} \left( E_k\right) \leq 4 x_k^{-1}
$
for all $k \geq k_1$. As $\sum_{k \in \mathbb{N}} x_k^{-1} < \infty$, the Borel-Cantelli lemma ensures that for almost all $\omega$, there exists $k_2(\omega)$ such that
$ \mathcal{P}_{x_k}  > \left \lfloor \phi_{k+1} \right \rfloor$ for all $k \geq k_2(\omega)$.

Applying Lemma \ref{lempoissomsc}, we deduce that, if we define the event $\Omega_0$ by 
$$
\Omega_0 := \left\{ \omega : \textrm{there exists $k_3(\omega)$ such that $N_{ \left \lfloor \phi_{k+1} \right \rfloor, j \phi_{k}^{\alpha}}(k) = 0$ for all $k \geq k_3(\omega)$} \right\},
$$
then $\mathbb{P}(\Omega_0) = 1$. Notice that there exists an integer $k_4$ such that for all $k \geq k_4$, 
$$
k^{-u'} \geq \exp \left ( - k \frac{1}{2(1-\alpha)} \frac{- \ln \operatorname*{L}(\theta)}{\theta} \right ).
$$ There exists a rank $k_5$ from which the sequence $(\phi_k)_{k \geq k_5}$ is increasing. Furthermore, that sequence tends to infinity. Let $\omega \in \Omega_0$. Let $n$ be an integer greater than $\phi_{k_5}$ large enough so that the unique $k \geq k_5$ satisfying $\phi_k \leq n < \phi_{k+1}$ is greater than $k_3(\omega)$ and $k_4$. Because $N_{ \left \lfloor \phi_{k+1}  \right \rfloor, j \phi_k^{\alpha}}(k) = 0$ and $n \leq \left \lfloor \phi_{k+1} \right \rfloor$, $N_{n, j \phi_k^{\alpha}}(k)=0$. Now, $j n^\alpha \geq j \phi_k^\alpha $, so  $N_{n,j n^\alpha}(k)=0$ and $H_{n,j n^\alpha} \leq k $. Moreover, as $n \geq \phi_k$, 
$$
k \leq (1 - \alpha) \frac{\theta}{- \ln \operatorname*L(\theta)} \ln n + u' (1-\alpha) \frac{\theta}{-\ln \operatorname*L(\theta)} \ln k.
$$
Further, as $k \geq k_4$, we have{~}:
$$
\exp \left (  k \frac{1}{2(1-\alpha)} \frac{- \ln \operatorname*L(\theta)}{\theta} \right )  \leq k^{-u'}  \exp \left (  k \frac{1}{1-\alpha} \frac{- \ln \operatorname*L(\theta)}{\theta} \right ) \leq n,
$$
so
$
k \leq  2 (1-\alpha) \frac{\theta}{-\ln \operatorname*L(\theta)} \ln n.
$
Therefore
$$
H_{n, j n^\alpha} \leq k \leq  (1 - \alpha) \frac{\theta}{-\ln \operatorname*L(\theta)} \ln n + u' (1-\alpha) \frac{\theta}{- \ln \operatorname*L(\theta)} \ln \left ( 2 (1-\alpha) \frac{\theta}{-\ln \operatorname*L(\theta)} \ln n \right ).
$$
We conclude that 
$$
\limsup_{n \rightarrow \infty} \frac{H_{n, j n^ \alpha} - (1 - \alpha) \frac{\theta}{- \ln \operatorname*L(\theta)} \ln n }{ \ln \ln n } \leq  u' (1-\alpha) \frac{\theta}{-\ln \operatorname*L(\theta)} \quad \textrm{a.s.}
$$
The fact that $u'$ can be chosen arbitrarily close to $\frac{1}{(1-\alpha) \theta}$ completes the proof.
\hfill $\square$

\subsection{Lower bound} 

We now study the lower bound. Combined with the upper bound supplied by Proposition~\ref{deuxxpromsc}, Theorem~\ref{coucoudidi} and Proposition~\ref{mscjetaime} will be proved.

\subsubsection{The case $j \geq \theta^{\ast}$ or $\alpha \in (0,1)$ and $\theta^{\ast} < \infty$} \label{sectdownone}

In this section, we prove that $H_{n,j} \geq C^{\ast} \ln n + O (\ln \ln n)$ a.s.\,\,if $j \geq \theta^{\ast}$ and (\ref{deuxhypohushi}) holds. We also prove that $H_{n, n^{\alpha}} \geq (1- \alpha) C^{\ast} \ln n + O (\ln \ln n)$ a.s.\,\,if $\theta^{\ast} < \infty$ and (\ref{deuxhypohushi}) holds. We shall see that the largest box plays a key role.

\begin{prop} \label{lowerno1poup} 
We suppose that $\theta^{\ast} < \infty$. Let $j \geq 1$ be an integer and $\alpha \in [0, 1)$ such that $(j, \alpha) \neq (1, 0)$. Under the assumption (\ref{deuxhypohushi}), we have~:
$$
\liminf_{n \rightarrow \infty} \frac{H_{n,j n^ \alpha} - (1 - \alpha) C^{\ast} \ln n }{ \ln \ln n } \geq \frac{3}{2 \ln \operatorname*{L}(\theta^{\ast})} \quad \textrm{a.s.}
$$
\end{prop}

\ \\
By definition, $H_{n, j n^{\alpha}} > k$ if and only if at generation $k$, there exists a box containing at least $j n^{\alpha}$ balls when $n$ balls have been thrown. We shall see that in our setting, it suffices to consider the largest box. As we intend to show that $H_{n, j n^{\alpha}}$ is bounded from below by $\left ( (1- \alpha) C^{\ast} - \varepsilon \right ) \ln n$, we take $k \approx \left ( (1- \alpha) C^{\ast} - \varepsilon \right ) \ln n$, \emph{i.e.}\,\,one initially throws $n \approx \exp \left ( k \left ( \frac{1}{(1- \alpha)C^{\ast}} + \varepsilon   \right ) \right )$ balls and we show that the largest box of the $k$-th generation contains at least $j n^{\alpha}$ balls. 

Let $\gamma$, $\gamma'$ and $\gamma''$ be three real numbers such that
$$
\gamma > \frac{3}{2 \theta^{\ast}}, \quad \gamma' > \frac{\gamma}{1 - \alpha} \quad \textrm{and} \quad \gamma' < \gamma'' < \frac{\gamma' - \gamma}{\alpha}.
$$
Note that if $\alpha = 0$, the third condition is simply $\gamma'' > \gamma'$. Define for all $k \in \mathbb{Z}_+$
$$
x_k := k^{\gamma'} \exp \left (  k \frac{1}{(1-\alpha) C^{\ast}}  \right ) \quad \textrm{and} \quad
\phi_k := k^{\gamma''} \exp \left  (  k \frac{1}{(1-\alpha) C^{\ast}}  \right ).
$$ 
We first show that a.s, $\mathcal{N}_{x_k, j \phi_k^{\alpha}}(k) \geq 1$ for all integers $k$ sufficiently large. To do so, we simply consider the largest box.

\begin{lemme} \label{lempoissoopppoki}
For almost all $\omega$, there exists $k_0(\omega)$ such that
$$
\mathcal{N}_{x_k, j \phi_k^{\alpha}}(k) \geq 1, \quad \textrm{for all $k \geq k_0(\omega)$}.
$$
\end{lemme}
\textbf{Proof :} 
For every generation $k \in \mathbb{N}$, we consider an imaginary box $\mathtt{b}(k)$ contained in the largest box $\overline{b}(k)$ of size $\overline{p}(k)$ such that a ball fallen in $\overline{b}(k)$ is thrown in $\mathtt{b}(k)$ with probability
$$
\frac{\mathtt{p}(k)}{\overline{p}(k)} \wedge 1, \quad \textrm{where} \quad \mathtt{p}(k) :=  k^{-\gamma} e^{-k / C^{\ast}}.
$$
Informally, the box $\mathtt{b}(k)$ has size $\mathtt{p}(k) \wedge \overline{p}(k)$. We denote by $A_k$ the event $ \{ \mathtt{p}(k) < \overline{p}(k) \}$ and by $B_k$ the event defined by 
$$
B_k := \{ \textrm{the box  $\mathtt{b}(k)$ contains less than $j \phi_k^\alpha$ balls when the first $\mathcal{P}_{x_k}$ have been thrown}    \}.
$$
Conditionally on $A_k$, the box $\mathtt{b}(k)$ has size $\mathtt{p}(k)$ so the number of balls contained in $\mathtt{b}(k)$ when $\mathcal{P}_{x_k}$ balls have been thrown is a Poisson variable with parameter $x_k \mathtt{p}(k)$. As a result
$$
\mathbb{P}(A_k \cap B_k) \leq  \mathbb{P}( B_k | A_k) = \mathbb{P} ( \mathcal{P}_{x_k \mathtt{p}(k)} < j \phi_k^{\alpha}  ).
$$
Applying (\ref{numbtwolemmod}), we get :
$$
\mathbb{P}(A_k \cap B_k) \leq d(p) j^{p} k^{-2},
$$
where   $p := 2/(\gamma' - \gamma - \alpha \gamma'') > 0$, and hence $\sum \mathbb{P} (A_k \cap B_k) < \infty$. By the Borel-Cantelli lemma, we deduce that a.s., for all $k$ sufficiently large,  $\mathtt{p}(k) \geq \overline{p}(k) $ or the box  $\mathtt{b}(k)$, which is contained in the box $\overline{b}(k)$, has at least $j \phi_k^\alpha$ balls when the first $\mathcal{P}_{x_k}$ balls have been thrown. Now, by Lemma \ref{thelimsuphushi}, we know that a.s., for all integers $k$ sufficiently large, $\mathtt{p}(k) < \overline{p}(k) $.
Lemma \ref{lempoissoopppoki} is therefore proved.  
\hfill $\square$

\ \\
We now deduce from Lemma \ref{lempoissoopppoki} the lower bound of $H_{n, j n^{\alpha}}$.

\ \\
\textbf{Proof of Proposition~\ref{lowerno1poup} :}
The same calculations performed at the beginning of the proof of Proposition~\ref{deuxxpromsc} show that for almost all $\omega$, there exists $k_1(\omega)$ such that
$\mathcal{P}_{x_k}  < \left \lceil \phi_{k-1} \right \rceil$ for all $k \geq k_1(\omega)$. Applying Lemma \ref{lempoissoopppoki}, we deduce that, if we define the event $\Omega_0$ by 
$$
\Omega_0 := \left\{ \omega : \textrm{there exists $k_2(\omega)$ such that $N_{ \left \lceil \phi_{k-1} \right \rceil, j \phi_k^{\alpha}}(k) \geq 1$ for all $k \geq k_2(\omega)$} \right\},
$$
then $\mathbb{P}(\Omega_0) = 1$. Let $\omega \in \Omega_0$. The sequence $(\phi_k)$ is increasing and tends to infinity. Let $n$ be an integer large enough so that the unique integer $k$ satisfying $\phi_{k-1} < n \leq \phi_k$ is greater than $k_2(\omega)$. Because $N_{ \left \lceil \phi_{k-1} \right \rceil, j \phi_{k}^{\alpha}}(k) \geq 1$ and $ \left \lceil \phi_{k-1} \right \rceil \leq n$, $N_{ n, j \phi_{k}^{\alpha}}(k) \geq 1$. Now, $j n^{\alpha} \leq j \phi_k^{\alpha}$, so $N_{ n, j n^{\alpha}}(k) \geq 1$ and $H_{n,j n^{\alpha}} > k.$ Further, as $n \leq \phi_k$,
$$
k \geq (1-\alpha)C^{\ast} \ln n - \gamma'' (1 - \alpha) C^{\ast} \ln k.
$$
As $\exp \left ( (k-1) \frac{1}{(1-\alpha)C^{\ast}} \right ) \leq \phi_{k-1} \leq n$, we have~:
$
k \leq (1-\alpha) C^{\ast} \ln n + 1.
$
Thus
$$
H_{n,j n^{\alpha}} > k \geq (1-\alpha)C^{\ast} \ln n - \gamma'' (1 - \alpha) C^{\ast} \ln \left (  (1-\alpha) C^{\ast} \ln n + 1  \right )
$$
and
$$
\liminf_{n \rightarrow \infty} \frac{H_{n,j n^{\alpha}}  -  (1-\alpha)C^{\ast} \ln n }{ \ln \ln n } \geq -  \gamma'' (1 - \alpha) C^{\ast} \quad \textrm{a.s.}
$$
The fact that $\gamma''$ can be chosen arbitrarily close to $\frac{3}{2 (1-\alpha) \theta^{\ast}}$ completes the proof.
\hfill $\square$

\subsubsection{The case $2 \leq j < \theta^{\ast}$ and $\alpha = 0$} \label{sectdowntwo}

We now prove that $H_{n,j} \geq C_j \ln n + O(\ln \ln n)$ a.s.\,\,if $j < \theta^{\ast}$. We shall notice that contrary to the case $j \geq \theta^{\ast}$, the number of boxes matters more than their sizes.

\begin{prop} \label{propberrou}
Let $j < \theta^{\ast}$ be an integer greater than 1. Then
$$
\liminf_{n \rightarrow \infty} \frac{H_{n,j} - \frac{j}{- \ln \operatorname*{L}(j)} \ln n }{ \ln \ln n } \geq  \frac{1}{2 \ln \operatorname*{L}(j)} \quad \textrm{a.s.}
$$
\end{prop}

\ \\
The main difference with the case $j \geq \theta^{\ast}$ is that we cannot consider the largest box any more. If $\theta^{\ast} < \infty$, one can even show by performing the usual calculations that with probability one, for all $k \geq C^{\ast} \ln n + o (\ln n)$, the largest box of the $k$-th generation contains no ball when $n$ balls have been initially thrown. We shall rather focus on some other boxes which are smaller but sufficiently numerous so that it is very unlikely that all of them contain less than $j$ balls when $n$ balls have been thrown. As we want to prove that $H_{n,j}$ is bounded from below by $\left ( \frac{j}{- \ln \operatorname*{L}(j)} - \varepsilon\right ) \ln n$, we consider the situation at the $k$-th generation when one initially throws approximatively $\exp \left ( k \left ( - \ln \operatorname*{L}(j) / j + \varepsilon \right ) \right )$ balls.

The boxes that will play a key role are those appearing in Lemma \ref{martibigg2} (recall that $j < \theta^{\ast}$){~}: with probability one,
$$
\lim_{k \rightarrow \infty}
\sqrt{k} e^{-k \varphi(j)} \# \left \{ i \in \mathbb{N}^k :  \mathtt{s}(k)   \leq p_i \leq 2\mathtt{s}(k)  \right \} = Q(j),
$$
where
$$
\mathtt{s}(k) :=  \exp \left ( - k \frac{ - {\operatorname*{L}}' (j)}{ \operatorname*{L}(j)} \right )
$$
and
$$
Q(j) := \frac{1}{\sqrt{2 \pi}} \frac{1 - 2^{-j}}{j} \left ( \frac{\operatorname*{L}''(j)}{\operatorname*{L}(j)}  - \left (  \frac{\operatorname*{L}'(j)}{\operatorname*{L}(j)} \right )^2  \right )^{-1/2} W(j) .
$$
Notice that, by Lemma \ref{martibigg}, as $j < \theta^{\ast}$, $Q(j) > 0$ a.s. Define
$$
\nu(k) :=  \left \lceil Q(j) k^{-1/2} e^{k \varphi(j)} / 2 \right \rceil.
$$
Then (recall that $\varphi(j) > 0$), for almost all $\omega$, there exists $k_0(\omega)$ such that
\begin{equation} \label{bertroua}
 \# \left \{ i \in \mathbb{N}^k : p_i \geq \mathtt{s}(k)  \right \} \geq \nu(k) \geq 1, \quad \textrm{for all $k \geq k_0(\omega)$}.
\end{equation}

We can now prove the following result.
\begin{lemme} \label{dernilemm}
Define for all $k \in \mathbb{Z}_+$
$$
x_k := k^u \exp \left ( k \frac{- \ln \operatorname*{L}(j)}{j} \right ), \quad \textrm{where} \quad u > \frac{1}{2 j}.
$$
Then for almost all $\omega$, there exists $k_1(\omega)$ such that
$$
\mathcal{N}_{x_k,j} (k) \geq 1, \quad \textrm{for all $k \geq k_1(\omega)$}.
$$
\end{lemme}
\textbf{Proof :} As $\mathbb{N}^{k}$ is countable, we can order the boxes of the $k$-th generation in the decreasing order. Denote by $\mathcal{G}_k$ the family of boxes having a size at least $\mathtt{s}(k)$. We create $\nu(k)$ new boxes{~}: if $\# \mathcal{G}_k \geq \nu(k)$, consider the first $\nu(k)$ boxes belonging to the family $\mathcal{G}_k$, denoted by $b_1(k), \dots, b_{\nu(k)} (k)$. For all $l \leq \nu(k)$, we place an imaginary box $\mathtt{b}_l(k)$ inside the box $b_l(k)$ such that a ball fallen in $b_l(k)$ is thrown in $\mathtt{b}_l(k)$ with probability $\mathtt{s}(k) / s_l(k)$, where $s_l(k)$ is the size of $b_l(k)$. In particular, every imaginary box has size $\mathtt{s}(k)$. If $\# \mathcal{G}_k < \nu(k)$,  we denote by $\mathtt{b}_1(k), \dots, \mathtt{b}_{\nu(k)} (k)$ the first $\nu(k)$ boxes. Introduce the events $A_k := \left \{ \# \mathcal{G}_k \geq \nu(k) \right \}$ and
$$
B_k := \left \{  \textrm{$\forall l \leq \nu(k)$, $\mathtt{b}_l(k)$ contains less than $j$ balls when the first $\mathcal{P}_{x_k}$ have been thrown}   \right \}.
$$
Conditionally on $\mathcal{F}_{\infty}$ and $A_k$, the boxes $\mathtt{b}_l(k)$ have size $\mathtt{s}(k)$ so the numbers of balls contained in $\mathtt{b}_l(k)$, $l \leq \nu(k)$, when $\mathcal{P}_{x_k}$ balls have been thrown are independent Poisson variables with parameters $x_k \mathtt{s}(k)$. Therefore
$$
\mathbb{P}(A_k \cap B_k | \mathcal{F}_{\infty}) \leq \mathbb{P}(\mathcal{P}_{x_k \mathtt{s}(k)} < j)^{\nu (k)}.
$$
Thus
$$
\ln \mathbb{P}(A_k \cap B_k | \mathcal{F}_{\infty}) \leq \frac{Q(j)}{2} k^{-1/2} e^{k \varphi(j)} \ln \mathbb{P}(\mathcal{P}_{x_k \mathtt{s}(k)} < j).
$$
It can be easily seen that $x_k \mathtt{s}(k)$ tends to 0, and that, as a result, 
$$
\ln \mathbb{P}(\mathcal{P}_{x_k \mathtt{s}(k)} < j) \sim - x_k^j \mathtt{s}(k)^j/ j!.
$$
Thus there exists a constant $c > 0$ such that for all $k \in \mathbb{Z}_+$,
$$
\ln \mathbb{P}(\mathcal{P}_{x_k \mathtt{s}(k)} < j) \leq - 2 c x_k^j \mathtt{s}(k)^j.
$$
Finally, we get
$$
\mathbb{P}(A_k \cap B_k | \mathcal{F}_{\infty}) \leq \exp \left (-  c Q(j) k^{-1/2} e^{k \varphi(k)} x_k^j \mathtt{s}(k)^j   \right ) = \exp \left (-  c Q(j) k^{j u -1/2}   \right ).
$$
As $u > \frac{1}{2 j}$ and $Q(j) > 0$ a.s., we get $ \mathbb{E} \left [ \sum \mathbf{1}_{A_k \cap B_k}  | \mathcal{F}_\infty \right] < \infty$ a.s., so  $\sum \mathbf{1}_{A_k \cap B_k} < \infty$ a.s. Combined with (\ref{bertroua}), this proves that a.s., for all integers $k$ sufficiently large, there exists an imaginary box containing at least $j$ balls when $\mathcal{P}_{x_k}$ have been thrown. As every imaginary box is contained in a real box, Lemma \ref{dernilemm} is proved.
\hfill $\square$

\ \\
We now deduce from Lemma \ref{dernilemm} the lower bound of $H_{n, j}$.

\ \\
\textbf{Proof of Proposition~\ref{propberrou} :}
Let $u$ and $u'$ be two real numbers such that $\frac{1}{2 j} < u < u'$. Define $ x_k := k^u \exp \left ( - k  \ln \operatorname*{L}(j) / j \right )$ and $\phi_k := k^{u'} \exp \left ( - k \ln \operatorname*{L}(j) / j \right )$. One can show that for almost all $\omega$, there exists $k_2(\omega)$ such that $ \mathcal{P}_{x_k}  < \left \lceil \phi_{k-1} \right \rceil$ for all $k \geq k_2(\omega)$. Applying Lemma \ref{dernilemm}, we deduce that if we define the event $\Omega_0$ by 
$$
\Omega_0 := \left\{ \omega : \textrm{there exists $k_3(\omega)$ such that $N_{ \left \lceil \phi_{k-1} \right \rceil, j }(k) \geq 1$ for all $k \geq k_3(\omega)$} \right\},
$$
then $\mathbb{P}(\Omega_0) = 1$. The sequence $(\phi_k)$ is increasing and tends to infinity. Let $\omega \in \Omega_0$. Let $n$ be an integer large enough so that the unique integer $k$ satisfying $\phi_{k-1} < n \leq \phi_{k}$ is greater than $k_3(\omega)$. Because $ N_{ \left \lceil \phi_{k-1} \right \rceil, j } (k) \geq 1$ and $n \geq \left \lceil \phi_{k-1} \right \rceil$, $N_{n,j} (k) \geq 1$, so $H_{n,j} > k $. Moreover, as $n \leq \phi_k$, 
$$
k \geq \frac{j}{- \ln \operatorname*{L}(j)} \ln n - u'  \frac{j}{- \ln \operatorname*{L}(j)} \ln k.
$$
As $n \geq \phi_{k-1} \geq \exp \left ( - (k-1) \ln \operatorname*{L}(j) / j \right )$, we have~: $k \leq \frac{j}{- \ln \operatorname*{L}(j)} \ln n +1$. Therefore 
$$
H_{n, j} > k \geq  \frac{j}{- \ln \operatorname*{L}(j)} \ln n - u'  \frac{j}{- \ln \operatorname*{L}(j)} \ln \left (  \frac{j}{- \ln \operatorname*{L}(j)} \ln n +1 \right ).
$$
We conclude that 
$$
\liminf_{n \rightarrow \infty} \frac{H_{n, j} -  \frac{j}{- \ln \operatorname*{L}(j)} \ln n }{ \ln \ln(n) } \geq - u'  \frac{j}{- \ln \operatorname*{L}(j)} \quad \textrm{a.s.}
$$
The fact that $u'$ can be chosen arbitrarily close to $\frac{1}{2 j}$ completes the proof.
\hfill $\square$

\begin{remark} \rm
Theorem~\ref{coucoudidi} states that $H_{n,j} \sim C_j \ln n$ a.s. This asymptotic behaviour was proved in the case $\theta^{\ast} = \infty$ (we even showed that $H_{n,j} = C_j \ln n+ O(\ln \ln n)$ a.s.). If $\theta^{\ast} < \infty$, then $-\ln \overline{p}(k) / k$ tends to $1 / C^{\ast}$ a.s.\,\,(see Lemma \ref{martibigg3}). We deduce from an argument similar to that in Proposition~\ref{lowerno1poup} that $H_{n,j} \sim C_j \ln n$ a.s. Finally, Theorem~\ref{coucoudidi} and Proposition~\ref{mscjetaime} have been proved.
\end{remark}

\section{Study of the saturation levels} \label{satusect}

In this section, we prove Theorem~\ref{maisonneuve} and Proposition~\ref{laroseara}. We shall first assume that (\ref{troishypohushi}) holds. We shall then study the case $\theta_{\ast} = -\infty$ to complete the proof of Theorem~\ref{maisonneuve}.

\subsection{The case $- \infty < \theta_{\ast} < 0$ and $\varphi(\theta_{\ast}) = 0$}

Throughout this section, we suppose that (\ref{troishypohushi}) holds. In particular, $C_{\ast} = -\theta_{\ast} / \ln \operatorname*{L} (\theta_{\ast}) = - \operatorname*{L} (\theta_{\ast}) / \operatorname*{L}' (\theta_{\ast})$. We are first interested in the lower bound. We shall then focus on the upper bound. We are inspired by the techniques developed in the previous section.

\subsubsection{Lower bound}

We show the inequalities (\ref{theoprenonkiooo}) and (\ref{theoprenonkioooij}) when (\ref{troishypohushi}) holds. Equation (\ref{numbtwolemmod}) will be the key tool.

\begin{prop} \label{deuxxpro} Suppose that (\ref{troishypohushi}) holds. Let $j \geq 1$ be an integer and $\alpha \in [0,1)$. Then
$$
\liminf_{n \rightarrow \infty} \frac{G_{n,j n^\alpha} - (1-\alpha) C_{\ast}  \ln(n) }{ \ln \ln(n) } \geq  -\frac{1}{\ln \operatorname*{L} (\theta_{\ast})} \quad \textrm{a.s.}
$$
\end{prop}

\ \\
Recall that $G_{n,j n^{\alpha}} > k$ if and only if at generation $k$, every box contains at least $j n^{\alpha}$ balls when $n$ balls have been thrown. Because we want to show that $G_{n, j n^{\alpha}}$ is bounded from below by  $\left ( (1 - \alpha) C_{\ast} - \varepsilon \right ) \ln(n)$, we take $k \approx \left (  (1 - \alpha) C_{\ast} - \varepsilon \right ) \ln(n)$, \emph{i.e.}\,\,$n \approx \exp \left (  k \left (  \frac{1}{(1 - \alpha) C_{\ast}} + \varepsilon \right ) \right )$.

\ \\
\textbf{Proof :} Let $u$ and $u'$ be two real numbers such that
$$
u > - \frac{1}{(1-\alpha) \theta_{\ast}} \quad \textrm{and} \quad u < u' < \frac{1}{\alpha} \left(  u + \frac{1}{\theta_{\ast}}   \right ),
$$
where the second condition reduces to $u' > u$ for $\alpha = 0$. Define for all $k \in \mathbb{N}$ :
$$
x_k := k^{u} \exp \left (  k \frac{1}{(1-\alpha) C_{\ast}}  \right ) \quad \textrm{and} \quad \phi_k := k^{u'} \exp \left (  k \frac{1}{(1-\alpha) C_{\ast}}  \right ).
$$ 
Proceeding as usual, let us prove the following key result{~}: for almost all $\omega$, there exists $k_0(\omega)$ such that
\begin{equation} \label{lempoisso}
\mathcal{M}_{x_k,j \phi_k^{\alpha}} (k)= 0 ,\quad \textrm{for all $k \geq k_0(\omega)$}.
\end{equation}
Let $x$ and $y$ be two real numbers such that $y \geq j$ and $k \in \mathbb{N}$. We calculate $\mathbb{E}\left[ \mathcal{M}_{x,y}(k) |\mathcal{F}_k \right] $. We write
$$
\mathcal{M}_{x,y}(k) = \sum_{i : |i|=k} \textbf{1}_{\mathcal{C}(i ; x) < y} =  \sum_{i : |i|=k} \textbf{1}_{\mathcal{C}(i ; x) < \lceil y \rceil}$$
Now, conditionally on $\mathcal{F}_k$,  $ (\mathcal{C}(i ; x))_{|i|=k} $ are Poisson variables with parameters $x p_i$, so
$$
\mathbb{E}\left[ \mathcal{M}_{x,y}(k)|\mathcal{F}_k \right]  =  \sum_{i : |i|=k}  \mathbb{P} (\mathcal{P}_{ x p_i} < \lceil y \rceil).
$$
Applying (\ref{numbtwolemmod}), we get
$$
\mathbb{E}\left[ \mathcal{M}_{x,y}(k)|\mathcal{F}_k \right]   \leq d(-\theta_{\ast}) \lceil y \rceil^{- \theta_{\ast}}  \sum_{i : |i|=k}   (x p_i)^{\theta_{\ast}} .
$$
Now, as $y \geq j$,
$$
\lceil y \rceil \leq \frac{y + 1}{y} y \leq \frac{j+1}{j} y.
$$
In the notations of Lemma \ref{martibigg}, we have :
$$
\mathbb{E}\left[ \mathcal{M}_{x,y}(k)|\mathcal{F}_k \right]   \leq d' y^{-\theta_{\ast}} x^{\theta_{\ast}} \operatorname*{L}(\theta_{\ast})^k W^{(k)}(\theta_{\ast}),
$$
where $d' := d(-\theta_{\ast}) (1 + 1/j)^{- \theta_{\ast}} $. We finally get
$$
\mathbb{E} \left[ \mathcal{M}_{x,y}(k) \right]  \leq d' y^{-\theta_{\ast}} x^{\theta_{\ast}} \operatorname*{L}(\theta_{\ast})^k.
$$
For $x= x_k$ and $y = j \phi_k^{\alpha} \geq j$, we obtain for every integer $k \in \mathbb{N}$ :
$$
\mathbb{E} \left[ \mathcal{M}_{x_k,j \phi_k^{\alpha}}(k) \right]  \leq d' j^{-\theta_{\ast}} k^{\theta_{\ast}(u - u' \alpha)}.
$$
Now $\theta_{\ast}(u - u' \alpha) < -1$. As a consequence
$$
\mathbb{E} \left[  \sum_{k \in \mathbb{Z}_+}  \mathbf{1}_{\mathcal{M}_{x_k,j \phi_k^{\alpha}}(k)  \geq 1}  \right] \leq  \mathbb{E} \left[  \sum_{k \in \mathbb{Z}_+}  \mathcal{M}_{x_k,j \phi_k^{\alpha}}(k)  \right] < \infty. 
$$
In particular, there is an a.s.\,\,finite number of integers $k$ such that $\mathcal{M}_{x_k,j \phi_k^{\alpha}}(k)  \geq 1$, which proves (\ref{lempoisso}).

Performing the same calculations as at the end of the proof of Proposition~\ref{lowerno1poup}, one finally gets{~}:
$$
\liminf_{n \rightarrow \infty} \frac{G_{n,j n^{\alpha}} -  (1-\alpha) C_{\ast} \ln n}{\ln \ln n} \geq - u' (1-\alpha) C_{\ast} \quad \textrm{a.s.}
$$
The fact that $u'$ can be chose arbitrarily close to $-\frac{1}{(1-\alpha) \theta_{\ast}}$ completes the proof. 
\hfill $\square$

\subsubsection{Upper bound}

In this section, we are interested in the upper bound of the saturation levels $G_{n,j}$. We prove that the inequalities (\ref{theoprenonkiooo}) and (\ref{theoprenonkioooij}) are in fact equalities whenever (\ref{deuxhypohushi}) and (\ref{troishypohushi}) hold by studying the smallest box, regardless of the value of $j$; there is no phase transition. 

\begin{prop} \label{deuxxprommm} Suppose that (\ref{deuxhypohushi}) and (\ref{troishypohushi}) hold. Let $j \geq 1$ be an integer and $\alpha \in (0,1)$. Then
$$
\limsup_{n \rightarrow \infty} \frac{G_{n,j} -C_{\ast}  \ln n }{ \ln \ln n } \leq \frac{3}{2 \ln \operatorname*{L}(\theta_{\ast})} +  \frac{C_{\ast}}{j} \quad \textrm{a.s.}
$$
and
$$
\limsup_{n \rightarrow \infty} \frac{G_{n,j n^{\alpha}} - (1-\alpha) C_{\ast}  \ln n }{ \ln \ln n } \leq \frac{3}{2 \ln \operatorname*{L}(\theta_{\ast})} \quad \textrm{a.s.}
$$
\end{prop}

\ \\
In order to prove both inequalities simultaneously, we may suppose that $\alpha \in [0,1)$, the case $G_{n,j}$ corresponding to $G_{n, j n^{\alpha}}$ with $\alpha = 0$.

By definition, $G_{n, j n^{\alpha}} \leq k$ if and only if at generation $k$, there exists a box containing less than $j n^{\alpha}$ balls when $n$ balls have been thrown. We shall see that it suffices to consider the smallest box. As we intend to show that $G_{n, j n^{\alpha}}$ is bounded from above by $\left ( (1- \alpha) C_{\ast} + \varepsilon \right ) \ln n$, we take $k \approx \left ( (1- \alpha) C_{\ast} + \varepsilon \right ) \ln n$, \emph{i.e.}\,\,$n \approx \exp \left ( k \left ( \frac{1}{(1- \alpha)C_{\ast}} - \varepsilon   \right ) \right )$.

\ \\
\textbf{Proof :} Let $\gamma$, $\gamma'$ and $\gamma''$ be three real numbers such that : 
\begin{itemize}
\item $\gamma > - \frac{3}{2 \theta_{\ast}}$, $\gamma' > \gamma + \frac{1}{j}$ and $\gamma'' > \gamma'$ if $\alpha = 0$, 
\item $\gamma > - \frac{3}{2 \theta_{\ast}}$, $\gamma' > \frac{\gamma}{1-\alpha}$ and $\gamma' < \gamma'' < \frac{\gamma' - \gamma}{\alpha}$ if $\alpha >0$.
\end{itemize}
Define
$$
x_k :=  k^{- \gamma'}  \exp \left ( k \frac{1}{ (1 - \alpha) C_{\ast}} \right ) \quad \textrm{and} \quad \phi_k :=   k^{- \gamma''}  \exp \left ( k \frac{1}{ (1 - \alpha) C_{\ast}} \right ).
$$
Let us prove the following result{~}: for almost all $\omega$, there exists $k_0(\omega)$ such that
\begin{equation} \label{lempoissoopppo}
\mathcal{M}_{x_k,j \phi_k^{\alpha}} (k) \geq 1, \quad \textrm{for all $k \geq k_0(\omega)$}.
\end{equation} 
Define  $\mathtt{p}(k) :=  k^ \gamma e^{-k / C_{\ast}}$. We consider an imaginary box $\mathtt{b}(k)$ at generation $k$ such that
\begin{enumerate}
\item if $1 \leq \mathtt{p}(k)$ : every ball thrown is placed in $\mathtt{b}(k)$.
\item if $\underline{p}(k) < \mathtt{p}(k) < 1$, where $\underline{p}(k)$ is the size of the smallest box $\underline{b}(k)$ : every ball fallen in $\underline{b}(k)$ is also placed in the imaginary box and every other ball is placed in $\mathtt{b}(k)$ with probability $(\mathtt{p}(k) - \underline{p}(k) )/(1- \underline{p}(k))$. Hence, the imaginary box has size $\mathtt{p}(k)$.
\item if $\mathtt{p}(k) \leq \underline{p}(k)$ : every ball fallen in the smallest box is placed in the imaginary box with probability $\mathtt{p}(k) / \underline{p}(k)$, and no other ball is placed in $\mathtt{b}(k)$.
\end{enumerate}
The box $\mathtt{b}(k)$ has thus size $\mathtt{p}(k) \wedge 1$, and whenever $\mathtt{p}(k) > \underline{p}(k)$, it contains the smallest box. From Lemma \ref{thelimsuphushi}, we know that a.s., for all integers $k$ sufficiently large, $\mathtt{p}(k) > \underline{p}(k)$. To prove Lemma  \ref{lempoissoopppo}, all we have to do is therefore to prove that $\mathbb{P}(\limsup A_k) = 0$, where the event $A_k$ is defined by
$$
A_k := \{  \textrm{$\mathtt{b}(k)$ contains at least $j \phi_k^{\alpha}$ balls when the first $\mathcal{P}_{x_k}$ balls have been thrown}   \}.
$$
By the Borel-Cantelli lemma, it suffices to show that $\sum \mathbb{P}(A_k) < \infty$. Let $k$ be an integer sufficiently large so that $\mathtt{p}(k) < 1$. The imaginary box has then size $\mathtt{p}(k)$, so we have :
$$
\mathbb{P}(A_k) = \mathbb{P} \left (   \mathcal{P}_{x_k \mathtt{p}(k)}   \geq j \phi_k^{\alpha} \right )
$$
We would like to apply (\ref{numbonelemmod}). To do so, let $p := j$ if $\alpha = 0$ and $p := 2 /(\gamma' -\gamma - \alpha \gamma'') > 0$ if $\alpha >0$. Note that for all integers $k$ sufficiently large, $j \phi_k^{\alpha} \geq p$. Equation (\ref{numbonelemmod}) then ensures that :
$$
\mathbb{P}(A_k) \leq c(p) j^{-p} k^{p (\alpha \gamma'' + \gamma - \gamma')}. 
$$
As $p (\alpha \gamma'' + \gamma - \gamma') < -1$, $\sum \mathbb{P} (A_k) < \infty$, which proves (\ref{lempoissoopppo}). 

Performing the same calculations as at the end of the proof of Proposition~\ref{deuxxpromsc}, one finally gets{~}:
$$
\limsup_{n \rightarrow \infty} \frac{G_{n,j n^{\alpha}} - (1-\alpha) C_{\ast} \ln n }{ \ln \ln n } \leq \gamma'' (1-\alpha) C_{\ast} \quad \textrm{a.s.}
$$
The fact that $\gamma''$ can be chosen arbitrarily close either to $ -\frac{3}{2 \theta_{\ast}} + \frac{1}{j}$ if $\alpha = 0$ or to 
$-\frac{3}{2 (1-\alpha) \theta_{\ast}}$ if $\alpha > 0$ completes the proof.
\hfill $\square$

\ \\
Proposition~\ref{laroseara} and a part of Theorem~\ref{maisonneuve} have been proved. We now turn our attention to the other part of Theorem~\ref{maisonneuve}.

\subsection{The case $\theta_{\ast} = - \infty$}

In this section, we prove that if $\theta_{\ast} = - \infty$, then $G_{n,j} \sim C_{\ast} \ln n$ a.s. We begin by showing the following lemma :

\begin{lemme} \label{lemjosko} If $\theta_{\ast} = - \infty$, then $ - \theta /  \ln \operatorname*{L} ( \theta)$ tends to $C_{\ast}$ as $\theta$ tends to $-\infty$. 
\end{lemme}

\ \\
\textbf{Proof :} The condition $\theta_{\ast} = - \infty$ means that $\varphi(\theta) > 0$ for all $\theta < 0$. Let $\psi = \ln \operatorname*{L}$. Recall that $C_{\ast} = \lim_{\theta \rightarrow -\infty} -1/ \psi'(\theta)$.  As $\psi$ is convex decreasing, it is known that $- \psi(\theta) / \theta$ tends to $l \in (0, \infty ]$ as $\theta$ tends to $ - \infty$. We distinguish two cases.

Either $l$ is finite, then one can easily show that $\theta \mapsto \psi (\theta) + l \theta$ is increasing. As a result, $\psi ' \geq - l$. If $l \neq  1 / C_{\ast}$, then $\psi'(\theta) \nrightarrow -l$ as $\theta$ tends to $- \infty$ , so there exists $\varepsilon > 0$ such that $\psi ' \geq \varepsilon - l$. Let $\theta < 0$. Then $\int_{\theta}^0 \psi' \geq (\varepsilon - l) (- \theta)$, \emph{i.e.}\,\,$\psi (\theta) + l \theta \leq \psi(0) + \varepsilon \theta$. Dividing by $- \theta$ and taking the limit, we get $\varepsilon \leq 0$, which is absurd. We deduce that $l = 1/C_{\ast}$, which means that $- \theta /  \ln \operatorname*{L} ( \theta) \rightarrow C_{\ast}$ as $\theta$ tends to $- \infty$.

Or $l$ is infinite. As for all $\theta < 0$, $\varphi(\theta) > 0$, the function $\theta \in (-\infty, 0) \mapsto \theta / \psi(\theta) - 1 / \psi'(\theta)$ is increasing. It is also positive. If it does not tend to 0 at $- \infty$, then it is bounded from below by some $\varepsilon > 0$. Multiplying by $- \psi'(\theta) > 0$, we get : $1 \geq 1 - \theta \psi ' (\theta) / \psi(\theta) \geq - \varepsilon \psi'(\theta)$. Consequently $\psi ' (\theta) \geq -1/ \varepsilon$  and $\psi(0) - \psi(\theta) \geq \theta / \varepsilon$. Dividing by $\theta < 0$, we obtain : $- \psi(\theta) / \theta \leq 1 / \varepsilon - \psi(0)/ \theta$. Taking the limit, we have $\infty \leq 1/ \varepsilon$, which is absurd. Finally, $\theta / \psi(\theta) - 1/ \psi ' (\theta)$ tends to 0 as $\theta$ tends to $-\infty$. Now, by definition of $C_{\ast}$, $\psi ' (\theta)$ tends to $-1/ C_{\ast}$. As a result, $- \theta / \psi(\theta)$ tends to $C_{\ast}$ as $\theta$ tends to $- \infty$. 
\hfill $\square$

\subsubsection{Lower bound}

We show that $G_{n,j} \geq C_{\ast} \ln n + o (\ln n)$ a.s. Equation (\ref{numbtwolemmod}) will be very useful.

\begin{prop} \label{huityaso} Suppose that $\theta_{\ast} = - \infty$. Let $j \geq 1$ be an integer. Then
$$
\liminf_{n \rightarrow \infty} \frac{G_{n,j}}{  \ln n } \geq C_{\ast} \quad \textrm{a.s.}
$$
\end{prop}

\ \\
\textbf{Proof :} We follow the usual strategy. Let $\theta < 0$ and $u > - \ln \operatorname*{L}(\theta) / \theta$. Define for all $k \in \mathbb{Z}_+$, $x_k := e^{ku}$.  Applying (\ref{numbtwolemmod}), we can show that
$$
\mathbb{E} [ \mathcal{M}_{x_k, j} (k)  | \mathcal{F}_k] \leq d(- \theta) j^{-\theta} x_k^{\theta} \operatorname*{L}(\theta)^k W^{(k)}(\theta),
$$
so
$$
\mathbb{E} [ \mathcal{M}_{x_k, j} (k) ] \leq d(- \theta) j^{-\theta}  \exp \left \{ k \theta ( u + \ln  \operatorname*{L}(\theta) / \theta )  \right \}.
$$
As $ \theta ( u + \ln  \operatorname*{L}(\theta) / \theta ) < 0$, $\mathbb{E} \left [ \sum  \mathcal{M}_{x_k, j} (k) \right ] < \infty$, so  a.s., for all integers $k$ sufficiently large, $\mathcal{M}_{x_k, j} (k) = 0$. The usual calculations yield{~}:
$$
\liminf_{n \rightarrow \infty} \frac{G_{n,j}}{\ln n} \geq \frac{- \theta}{\ln \operatorname*{L}(\theta)} \quad \textrm{a.s.}
$$
Lemma \ref{lemjosko} enables us to complete the proof.  
\hfill $\square$

\subsubsection{Upper bound}

We finally prove that $G_{n,j} \leq C_{\ast} \ln n + o(\ln n)$ a.s. As $\theta_{\ast} = - \infty$, Lemma \ref{martibigg3} cannot be applied (we do not know how the size of the smallest box behaves). Surprisingly, we are inspired by the proof of Proposition~\ref{propberrou}.

\begin{prop} \label{huityasoik} Suppose that $\theta_{\ast} = - \infty$. Let $j \geq 1$ be an integer. Then
$$
\limsup_{n \rightarrow \infty} \frac{G_{n,j}}{  \ln n } \leq C_{\ast} \quad \textrm{a.s.}
$$
\end{prop}

\ \\
\textbf{Proof :} We may suppose that $j=1$.  Let $\theta < 0$. Here, the boxes of the $k$-th generation that will play a key role are those having a size approximatively $\mathtt{s}(k)$, where
$$
\mathtt{s}(k) :=  \exp \left ( - k \frac{ - {\operatorname*{L}}' (\theta)}{ \operatorname*{L}(\theta)} \right ).
$$
Applying Lemma \ref{martibigg2}, we have with probability one :
$$
\lim_{k \rightarrow \infty}
\sqrt{k} e^{-k \varphi(\theta)} \# \left \{ i \in \mathbb{N}^k :  \mathtt{s}(k) / 2   \leq p_i \leq \mathtt{s}(k)  \right \} = Q(\theta),
$$
where
$$
Q(\theta) := \frac{1}{\sqrt{2 \pi}} \frac{1 - 2^{\theta} }{- \theta} \left ( \frac{\operatorname*{L}''(\theta)}{\operatorname*{L}(\theta)}  - \left (  \frac{\operatorname*{L}'(\theta)}{\operatorname*{L}(\theta)} \right )^2  \right )^{-1/2} W(\theta) .
$$
Notice that, by Lemma \ref{martibigg}, $Q(\theta) > 0$ a.s. Define
$$
\nu(k) :=  \left \lceil Q(\theta) k^{-1/2} e^{k \varphi(\theta)} / 2 \right \rceil.
$$
Then (recall that $\varphi(\theta) > 0$), for almost all $\omega$, there exists $k_0(\omega)$ such that
\begin{equation} \label{bertroualol}
\# \left \{ i \in \mathbb{N}^k : p_i \leq \mathtt{s}(k)  \right \} \geq \nu(k) \geq 1, \quad \textrm{for all $k \geq k_0(\omega)$}.
\end{equation}

Let $u < -\operatorname*{L}'(\theta)/\operatorname*{L}(\theta)$. Define for all $k \in \mathbb{Z}_+$, $x_k := e^{k u}$. Let us prove that for almost all $\omega$, there exists $k_1(\omega)$ such that
\begin{equation} \label{dernresk}
\mathcal{M}_{x_k,1} (k) \geq 1, \quad \textrm{for all $k \geq k_1 (\omega)$}. 
\end{equation}
Denote by $\mathcal{G}_k$ the family of boxes having a size at most $\mathtt{s}(k)$. We define $\nu(k)$ new boxes~: if $\# \mathcal{G}_k \geq \nu(k)$, we create $\nu(k)$ imaginary boxes  $\mathtt{b}_1(k), \dots, \mathtt{b}_{\nu(k)} (k)$ of size $\mathtt{s}(k)$ containing the first $\nu(k)$ real boxes belonging to the family $\mathcal{G}_k$. If $\# \mathcal{G}_k < \nu(k)$, we denote by $\mathtt{b}_1(k), \dots, \mathtt{b}_{\nu(k)} (k)$ the first $\nu(k)$ boxes of the $k$-th generation. Introduce the event $A_k := \left \{ \# \mathcal{G}_k \geq \nu(k) \right \}$ and
$$
B_k := \left \{  \textrm{$\forall l \leq \nu(k)$, $\mathtt{b}_l(k)$ contains at least one ball when the first $\mathcal{P}_{x_k}$ have been thrown}   \right \}.
$$
Conditionally on $\mathcal{F}_{\infty}$ and $A_k$, the boxes $\mathtt{b}_l(k)$ have size $\mathtt{s}(k)$ so the numbers of balls contained in $\mathtt{b}_l(k)$, $l \leq \nu(k)$, when $\mathcal{P}_{x_k}$ balls have been thrown are independent Poisson variables with parameters $x_k \mathtt{s}(k)$. Therefore
$$
\mathbb{P}(A_k \cap B_k | \mathcal{F}_{\infty}) \leq \mathbb{P}(\mathcal{P}_{x_k \mathtt{s}(k)} \geq 1)^{\nu (k)}.
$$
Thus
$$
\ln \mathbb{P}(A_k \cap B_k | \mathcal{F}_{\infty}) \leq \frac{Q(\theta)}{2} k^{-1/2} e^{k \varphi(\theta)} \ln \mathbb{P}(\mathcal{P}_{x_k \mathtt{s}(k)} \geq 1).
$$
Now $\mathbb{P}(\mathcal{P}_{x_k \mathtt{s}(k)} \geq 1) \leq x_k \mathtt{s}(k)$, so
$$
 \mathbb{P}(A_k \cap B_k | \mathcal{F}_{\infty}) \leq \exp \left \{ \left ( u + \frac{\operatorname*{L}'(\theta)}{\operatorname*{L}(\theta)} \right )  \frac{Q(\theta)}{2} k^{1/2} e^{k \varphi(\theta)}  \right \}.
$$
As $u + \operatorname*{L}'(\theta) / \operatorname*{L}(\theta) < 0$, $\varphi(\theta) >0$  and $Q(\theta) > 0$ a.s., we get $ \mathbb{E} \left [ \sum \mathbf{1}_{A_k \cap B_k}  | \mathcal{F}_\infty \right] < \infty$ a.s., so  $\sum \mathbf{1}_{A_k \cap B_k} < \infty$ a.s. Combined with (\ref{bertroualol}), this proves that a.s., for all integers $k$ sufficiently large, there exists an empty imaginary box when $\mathcal{P}_{x_k}$ have been thrown. As every imaginary box contains a real box, (\ref{dernresk}) is proved.

Proceeding as usual, we deduce that for all $\theta < 0$,
$$
\limsup_{n \rightarrow \infty} \frac{G_{n,1}}{  \ln n } \leq \frac{\operatorname*{L}(\theta)}{- \operatorname*{L}'(\theta)} \quad \textrm{a.s.}
$$
We get the result by letting $\theta$ tend to $- \infty$. 
\hfill $\square$

\ \\
Finally, Theorem~\ref{maisonneuve} has been proved.

\section{An explanation for the phase transition} \label{studexap}

In this section, we explain heuristically why a phase transition may occur in the asymptotics of $H_{n,j}$ but not in those of $H_{n,n^{\alpha}}$, $G_{n,j}$ and $G_{n,n^{\alpha}}$. To do so, we rephrase the setting in terms of interval-fragmentations. One can easily construct a family $(F(k), k \in \mathbb{Z}_+)$ of random open subsets of $(0,1)$ with $F(0) = (0,1)$, which is nested in the sense that $F(k') \subseteq F(k)$ for all $k \leq k'$, and such that the set of the lengths of the interval components of $F(k)$ is $\{ p_i : |i| = k \}$ for every integer $k$. The boxes correspond to the interval components and the balls to a sequence $(U_i)_{i \in \mathbb{N}}$ of independent random variables uniformly distributed on $(0,1)$ and independent of the fragmentation $F$. The height $H_{n,j}$ corresponds to the least integer $k$ such that every interval component of $F(k)$ contains less than $j$ elements of the set $\{U_1, \dots, U_n\}$, and the saturation level $G_{n,j}$ to the least integer $k$ such that there exists an interval component of $F(k)$ containing less than $j$ elements of the set $\{U_1, \dots, U_n\}$. Roughly speaking, the height $H_{n,j}$ depends crucially on the minimal length $\underline{m}_{n,j}$ of the intervals $[\hat{U}_i, \hat{U}_{i+j-1}]$ for $1 \leq i \leq n-j+1$, where $0 < \hat{U}_1 < \dots < \hat{U}_n <1$ are the ordered statistics of the family $(U_1, \dots, U_n)$, whereas the saturation level $G_{n,j}$ is related to the maximal length $\overline{m}_{n,j}$ of the intervals $[\hat{U}_i, \hat{U}_{i+j}]$ for $0 \leq i \leq n-j+1$, where $\hat{U}_0 := 0$ and $\hat{U}_{n+1} :=1$. Indeed, the height $H_{n,j}$ is the first time $k$ when no cluster $[\hat{U}_{i}, \hat{U}_{i+j-1}]$, $1 \leq i \leq n-j+1$, of size $j$ (in particular the smallest one) is included in an interval component of $F(k)$, and the saturation level $G_{n,j}$ is the first time $k$ when there exists a cluster $[\hat{U}_{i}, \hat{U}_{i+j}]$, $0 \leq i \leq n-j+1$, of size $j+1$ (possibly the largest one) containing an interval component of $F(k)$.

It is easy to show from Lemma \ref{martibigg3} that if we used equidistributed points $\{ j/(n+1) : 1\leq j \leq n \}$ instead of i.i.d.\,\,uniform points, then no phase transition would occur. More precisely, all the heights $H_{n,j}$ would be equivalent to $C^{\ast} \ln n$ a.s. Further, the heights $H_{n, n^{\alpha}}$ would be equivalent to $(1-\alpha) C^{\ast} \ln n$ a.s.\,\,and the saturation levels $G_{n, j n^{\alpha}}$ to $(1-\alpha) C_{\ast} \ln n$ a.s.

We first explain why the clusters of size $n^{\alpha} $ behave as if the points $U_i$, $1 \leq i \leq n$, were equidistributed on $(0,1)$. It will follow that no phase transition occurs in the asymptotics of $H_{n, n^{\alpha}}$ and $G_{n, n^{\alpha}}$. To do so, let us first prove that $\underline{m}_{n,n^{\alpha}} \geq n^{-1+\alpha + o(1)}$ a.s.  Let $\varepsilon > 0$. We have~:
$$
\mathbb{P}(\underline{m}_{n,n^{\alpha}} \leq n^{-1+\alpha - 2 \varepsilon}) \leq n \mathbb{P}(\hat{U}_{n^\alpha} - \hat{U}_1 \leq n^{-1+\alpha - 2 \varepsilon}).
$$
Let $(e_i)_{i \in \mathbb{N}}$ be a sequence of independent exponential variables with parameters 1. Define $\gamma_k := \sum_{1 \leq i \leq k} e_i$, for all $k \geq 1$. Then we have~:
\begin{equation} \label{sidneyros}
(\hat{U}_1, \dots, \hat{U}_n) \overset{\textrm{(d)}}{=} (\gamma_1/\gamma_{n+1}, \dots, \gamma_n / \gamma_{n+1}),
\end{equation}
so 
\begin{eqnarray*}
\mathbb{P}(\hat{U}_{n^\alpha} - \hat{U}_1 \leq n^{-1+\alpha - 2 \varepsilon}) & = & \mathbb{P} \left ( \frac{\gamma_{n^{\alpha}}- \gamma_1}{\gamma_{n+1}}  \leq  n^{-1+\alpha - 2 \varepsilon} \right ) \\
& \leq & \mathbb{P} \left ( \gamma_{n^{\alpha}}- \gamma_1  <  n^{1+\varepsilon} n^{-1+\alpha - 2 \varepsilon} \right ) + \mathbb{P} \left ( \gamma_{n+1} \geq  n^{1+\varepsilon} \right ) \\
& = &  \mathbb{P} \left ( \mathcal{P}_{n^{\alpha - \varepsilon}} \geq n^{\alpha} \right ) + \mathbb{P} \left ( \mathcal{P}_{n^{1+\varepsilon}} \leq n+1 \right ) 
\end{eqnarray*}
 since the sequence $(\gamma_i)_{i \in \mathbb{N}}$ has the same law as the ordered sequence of points of a standard Poisson process. Applying Lemma \ref{lemmemodif} with $p = 3/\varepsilon$, we have for every integer $n$ sufficiently large~:
$$
\mathbb{P} \left ( \mathcal{P}_{n^{\alpha - \varepsilon}} \geq n^{\alpha} \right ) + \mathbb{P} \left ( \mathcal{P}_{n^{1+\varepsilon}} \leq n+1 \right ) \leq \left ( c(p) + 2^p d(p) \right ) n^{-3}.
$$
We deduce that that for every integer $n$ sufficiently large, we have~:
$$
\mathbb{P}(\underline{m}_{n,n^{\alpha}} \leq n^{-1+\alpha - 2 \varepsilon}) \leq  \left ( c(p) + 2^p d(p) \right ) n^{-2}.
$$
Applying  the Borel-Cantelli lemma and letting $\varepsilon$ tend to 0, we easily conclude that $\underline{m}_{n,n^{\alpha}} \geq n^{-1+\alpha + o(1)}$ a.s. Similar calculations yield~: $\overline{m}_{n,n^{\alpha}} \leq n^{-1+\alpha + o(1)}$ a.s. We conclude that the clusters of size $n^{\alpha}$ essentially behave as if the points $U_i$, $1 \leq i \leq n$, were equidistributed on $(0,1)$, which explains why no phase transition occur in the asymptotics of $H_{n, n^{\alpha}}$ and $G_{n, n^{\alpha}}$.

On the contrary, the clusters of size $j$ behave very differently from the equidistribution case. For instance, thanks to the identification (\ref{sidneyros}), we may apply well-known results on extreme values (see \emph{e.g.}\,\,\cite{resnick}). One finds that $\underline{m}_{n,j} \approx n^{-j/(j-1)}$, so that $\underline{m}_{n,j}$ is much smaller than $j/n$, even at a logarithmic scale. Furthermore, the larger the integer $j$ is, the closer these two quantities are, so that $H_{n,j}$ may eventually be equivalent to $C^{\ast} \ln n$ a.s.\,\,if $j$ is sufficiently large. On the other hand, if the integer $j$ is too small, then $H_{n,j}$ may be greater than $(1+\varepsilon) C^{\ast} \ln n + o(\ln n)$ a.s.\,\,for some $\varepsilon > 0$~: one has to wait for a long time before the smallest clusters of size $j$ are no longer included in any interval component of the fragmentation; in that case, there is a phase transition. 

Concerning the saturation level $G_{n,j}$, one gets $\overline{m}_{n,j} \approx \ln n / n$, which differs from the equidistribution case only up to a logarithmic factor, and the impact of the latter is asymptotically negligible. It explains why no phase transition occurs in the asymptotics of $G_{n,j}$.

\ \\
\textbf{Acknowledgements.} I would like to express my gratitude to my thesis director Jean Bertoin who suggested this subject to me. I sincerely thank him for his invaluable advice and for his careful reading.

\nocite{*}

\bibliographystyle{plain}
\bibliography{bibliarticle}

\end{document}